\documentclass[preprint,11pt]{elsarticle}

\usepackage{graphicx}
\usepackage{amssymb,amsmath}
\usepackage{subfigure}
\usepackage{float}

\newtheorem{lemma}{Lemma}
\newtheorem{rem}{Remark}
\newtheorem{definition}{Definition}
\newtheorem{example}{Example}[subsection]
\newproof{proof}{Proof}
\numberwithin{equation}{section}

\usepackage{lineno}
\usepackage{xcolor}
\newcommand{\vv}{\overrightarrow{V}}

\newcommand{\uAve}{\bar{u}}

\newcommand{\vW}{\textbf{W}}

\begin{document}

\begin{frontmatter}


\title{Cell-average based neural network method for hyperbolic and parabolic partial differential equations}



\author[label1]{Changxin Qiu}
\ead{cxqiu@iastate.edu}

\author[label1]{Jue Yan\corref{cor1}\fnref{label2}}
\ead{jyan@iastate.edu}
\fntext[label2]{Research work of the author is supported by National Science Foundation grant DMS-1620335 and Simons Foundation grant 637716.}
\cortext[cor1]{corresponding author}
\address[label1]{Department of Mathematics, Iowa State University, Ames, IA 50011, USA}

\begin{abstract}
Motivated by finite volume scheme, a cell-average based neural network method is proposed. The method is based on the integral or weak formulation of partial differential equations. A simple feed forward network is forced to learn the solution average evolution between two neighboring time steps. Offline supervised training is carried out to obtain the optimal network parameter set, which uniquely identifies one finite volume like neural network method. Once well trained, the network method is implemented as a finite volume scheme, thus is mesh dependent. Different to traditional numerical methods, our method can be relieved from the explicit scheme CFL restriction and can adapt to any time step size for solution evolution. For Heat equation, first order of convergence is observed and the errors are related to the spatial mesh size but are observed independent of the mesh size in time. The cell-average based neural network method can sharply evolve contact discontinuity with almost zero numerical diffusion introduced. Shock and rarefaction waves are well captured for nonlinear hyperbolic conservation laws. 

\end{abstract}


\end{frontmatter}

%\linenumbers

\section{Introduction}

In this paper, we develop cell-average based neural network (CANN) method solving time dependent hyperbolic and parabolic partial differential equations (PDEs)
\begin{equation}\label{intro:conv-diff}
   u_t+f(u)_x=\mu u_{xx}. 
\end{equation}
Our general idea is to follow available numerical schemes to build up neural network methods. We consider to combine the powerful machine learning mechanism of neural networks and the principles of classical numerical methods to develop suitable neural network methods for partial differential equations (\ref{intro:conv-diff}). 

Machine learning with neural networks have achieved tremendous success in image classification, text, videos and speech recognition \cite{LeCun-Bengio-1995, Bengio-2009, Krizhevsky-Sutskever-Hinton-2012, LeCun-Bengio-Hinton-2015} for the last three decades. In the last few years, connections between differential equations and machine learning have been established, i.e. \cite{E-2017, Chaudharip-Osher2017,Rudy-Brunton-Proctor-Kutz2017,  Chang-Meng-Haber-Ruthotto-Begert-Holtham-2018, Long-Lu-Dong-2019,  Ruthotto-Haber2020, He-Xu2019}. Very recently, machine learning neural networks have also been explored to directly solve partial differential equations, for which a better solver may be obtained or it can assist to improve the performance of current numerical methods.

Nonlinear convection diffusion equation (\ref{intro:conv-diff}) may not be complex and hard to solve in general. It can be used as a prototype or model equation for the more complicated Euler and Navier-Stokes equations for fluid dynamics. There exist quite a few numerical methods successfully developed for (\ref{intro:conv-diff}). Several numerical challenges are still present. For example, a lot of ongoing efforts are toward investigating implicit or semi-implicit methods to obtain efficient solvers. It turns out, once well trained, our cell-average or finite volume based neural network method can be relieved from the explicit scheme CFL restriction and can adapt large time step size (i.e. $\Delta t=4\Delta x$) for solution evolution, even being implemented as an explicit method.

\subsection{Related works}{\label{S:1.1}}

Depending on the ways of neural networks applied and the goals, neural network methods can be roughly classified into two groups. One group is to design and apply neural network methods to directly solve partial differential equations. The second group is to have neural networks applied to assist and improve available numerical methods. 

One popular class is to find best solution representation in terms of neural networks. With network input vector as $x$ and $t$, such methods have the advantage of automatic differentiation, mesh free and can be applied to solve many types of PDEs. Approximation power of neural networks \cite{Cybenko-1989,HORNIK-1990,Barron-1993,Leshno-Lin-1993,Pinkus-1999} are explored with these methods. For solving PDEs, we have the early works of \cite{Lagaris-Likas-Fortiadis1998,RUDD2015}, the works of \cite{Berg-Nystrom-2018, Sirignano-Spiliopoulos-2018} and the popular physics-informed (PINN) methods of \cite{Raissi-Perdikaris-Karniadakis-2017a,Raissi-Perdikaris-Karniadakis-2017b,Raissi-Perdikaris-Karniadakis-2019,Lu-Meng-Mao-Karniadakis-2021}. To improve the PINN efficiency on larger domain, an extreme and distributed network is considered in  \cite{Dwivedi-Srinivasan-2020}. Application to incompressible Navier-Stokes equations is studied in \cite{Jin-Cai-Li-Karniadakis-2021}. We also have the works of  \cite{Zang-Bao-Ye-Zhou-2020} and \cite{Cai-2021-1,Cai-2021-2}, for which weak formulations are applied in the loss function in stead of PDEs explicitly enforced on collocation points. Comparison to reduced basis method is studied in \cite{Kutyniok-2021} and performance comparison to finite volume and discontinuous Galerkin methods are considered in \cite{Michoski-Oliver-Hatch-2020}. We refer to \cite{Shin-Darbon-Karniadakis-2020,Laakmann-Petersen2021} for PINN convergence studies.

Neural network methods have been found rather successful for solving high dimensional PDEs. We refer to the early work of \cite{Beck-E-Jentzen-2019} and other discussions in \cite{Chan-Wai-Nam-2019,Pham-Warin-Germain-2021,Hutzenthaler-2020,Lye-Mishra-Ray2020}. Other studies include designing specific networks for elliptic type PDEs or solving inverse problems, see \cite{Fan-Lin-Ying-Zepeda-2019,khoo-lu-ying-2020,Li-Lu-Mao-2020}. In \cite{Winovich-Ramani-Lin2019}, convolutional networks are used for estimating the mean and variance for uncertainty quantification. We also have the works of \cite{Wu-Xiu-2020} and \cite{Qin-Chen-Jakeman-Xiu-2021}, for which method of lines approach is explored with Fourier basis considered and Residual networks applied to evolve the dynamical system. Different to PINN or related methods in which a global solution with neural networks is sought, our method is similar to finite volume scheme thus is mesh dependent and is a local solver.

The other group is to combine networks with classical numerical methods for performance improvement. We have early result of \cite{Ray-Hesthaven-2018} applying neural networks as trouble-cell indicator. Deep reinforcement network in \cite{Wang-Shen-Long-Dong-2019} is explored to estimate the weights and enhance WENO schemes performance. In \cite{Discacciati-Hesthaven-Ray-2020}, network is applied for identifying suitable amount of artificial viscosity added. Neural network of \cite{Hsieh-Zhao-2019} is applied to speed up iterative solver for elliptic type PDEs. Furthermore, we have WENO schemes augmented with convolutional networks for shock detection in \cite{Sun-Wang-Chang-Xing-Xiu-2020} and DG methods with imaging edge detection technique with convolutional network explored for strong shock detection in \cite{Beck-Zeifang-Schwarz-Flad2020}.
We further have the work of \cite{Yu-Shu-2021} in which network is applied estimating total variation bounded constants for DG methods.

\subsection{Motivation and our approach}{\label{S:1.2}}

The cell-average based neural network method is closely related to finite volume scheme. Let's review first order upwind finite volume method for linear advection equation
\begin{equation}\label{intro:1D-example}
 u_t+u_x=0,
\end{equation}
which motivates the design of our cell-average based method. Integrating the advection equation (\ref{intro:1D-example}) over one computational cell $I_j=[x_{j-1/2},x_{j+1/2}]$, with cell average notation $\bar{u}_j=\frac{1}{\Delta x}\int_{I_j} u(x,t) ~dx$ and numerical flux $\hat{u}$ introduced at $x_{j\pm 1/2}$, we obtain $\frac{\partial \bar{u}_j}{\partial t}+ \frac{\hat{u}_{j+1/2}-\hat{u}_{j-1/2}}{\Delta x}=0$. This is the starting point of finite volume schemes. Adapting upwind for the numerical flux and forward Euler for time discretization, we have the following first order upwind finite volume method 
\begin{equation}\label{intro:upwind-forward}
    \frac{\bar{u}^{n+1}_{j}-\bar{u}^{n}_{j}}{\Delta t}+\frac{\bar{u}^{n}_{j}-\bar{u}^{n}_{j-1}}{\Delta x}=0.
\end{equation}

For cell-average based neural network method, we consider to integrate the advection equation (\ref{intro:1D-example}) over rectangle box of $I_j\times(t_n,t_{n+1})$, which covers both the computational cell in space and the sub interval in time. The integral or weak formulation of the advection equation is obtained as 
\begin{equation}\label{intro:1D-integralFormat}
    \bar{u}_j(t_{n+1}) =\bar{u}_j(t_n)-\frac{1}{\Delta x}\int^{t_{n+1}}_{t_n}\int_{I_j} u_x~ dxdt.
\end{equation}
Exact solution of (\ref{intro:1D-example}) satisfies the above integral formulation. With extra regularity on the solution, the above integral format is equivalent to the advection equation (\ref{intro:1D-example}) given in differentiation format.
A simple feed forward neural network is adapted to approximate 
spatial variable related term 
$\mathcal{N}(\vv_{j}^{in};\Theta)\approx\int^{t_{n+1}}_{t_n}\int_{I_j} u_x ~dxdt$. Here $\vv_{j}^{in}$ is the network input vector and $\Theta$ is the network parameter set. And we obtain the following cell-average based neural network method for (\ref{intro:1D-example})
\begin{equation}\label{intro:CANN-scheme}
     \bar{v}_j^{n+1} = \bar{v}_j^{n}+\mathcal{N}(\vv_{j}^{in}; \Theta^{\star}), \quad \forall j\quad\mbox{and}\quad \forall n.
\end{equation}
Setting $\bar{v}_j^{n}=\bar{u}_j(t_n)$ and by minimizing the difference between network output $\bar{v}_j^{n+1}$ and the target $\bar{u}_j(t_{n+1})$, optimal parameter set $\Theta^{\star}$ can be obtained through offline supervised learning.

\begin{figure}[!htbp]
\centering
\subfigure[]{
\includegraphics[width=0.24\linewidth]{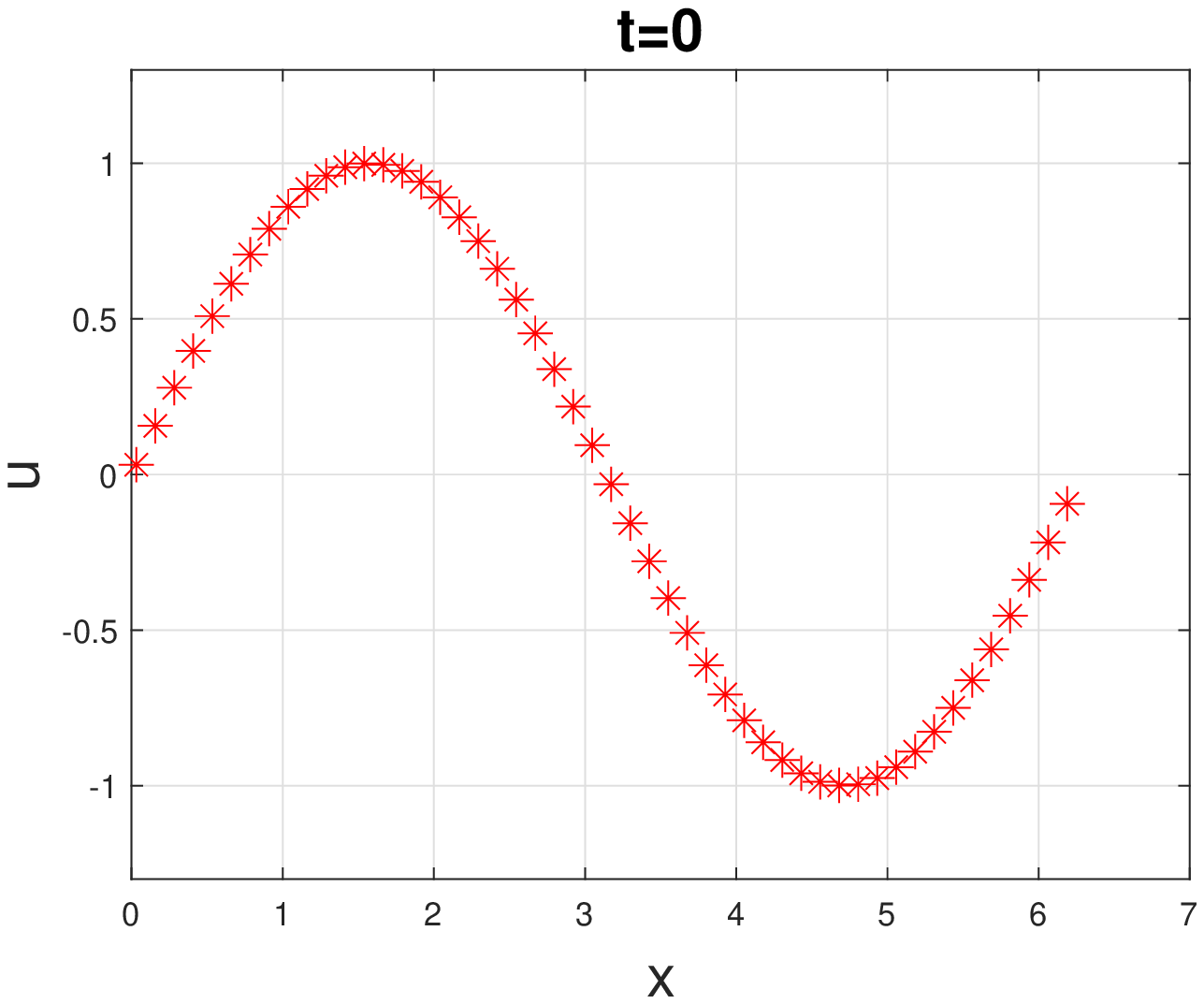}}\,
\subfigure[]{\includegraphics[width=0.24\linewidth]{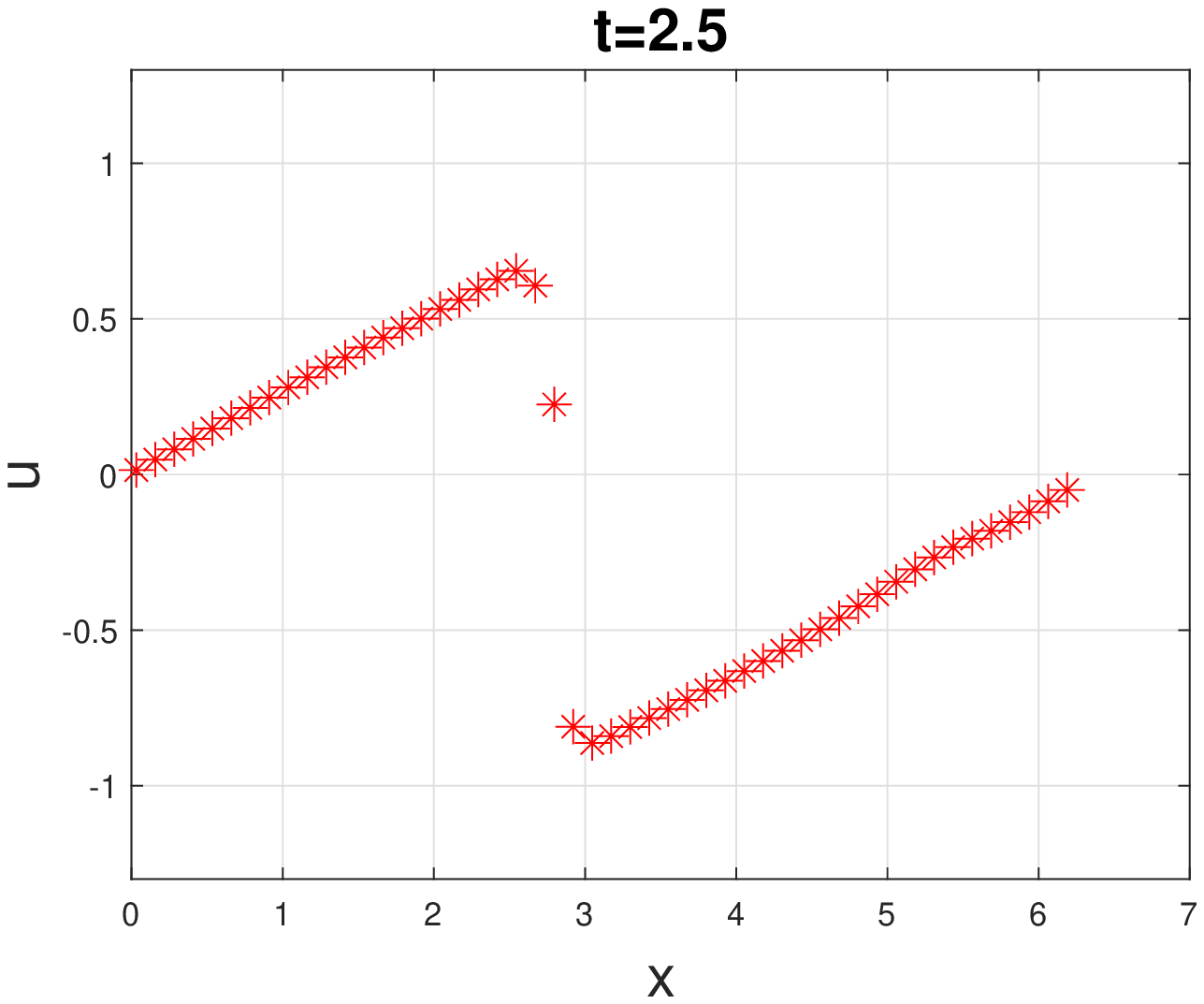}}\,
\subfigure[]{\includegraphics[width=0.24\linewidth]{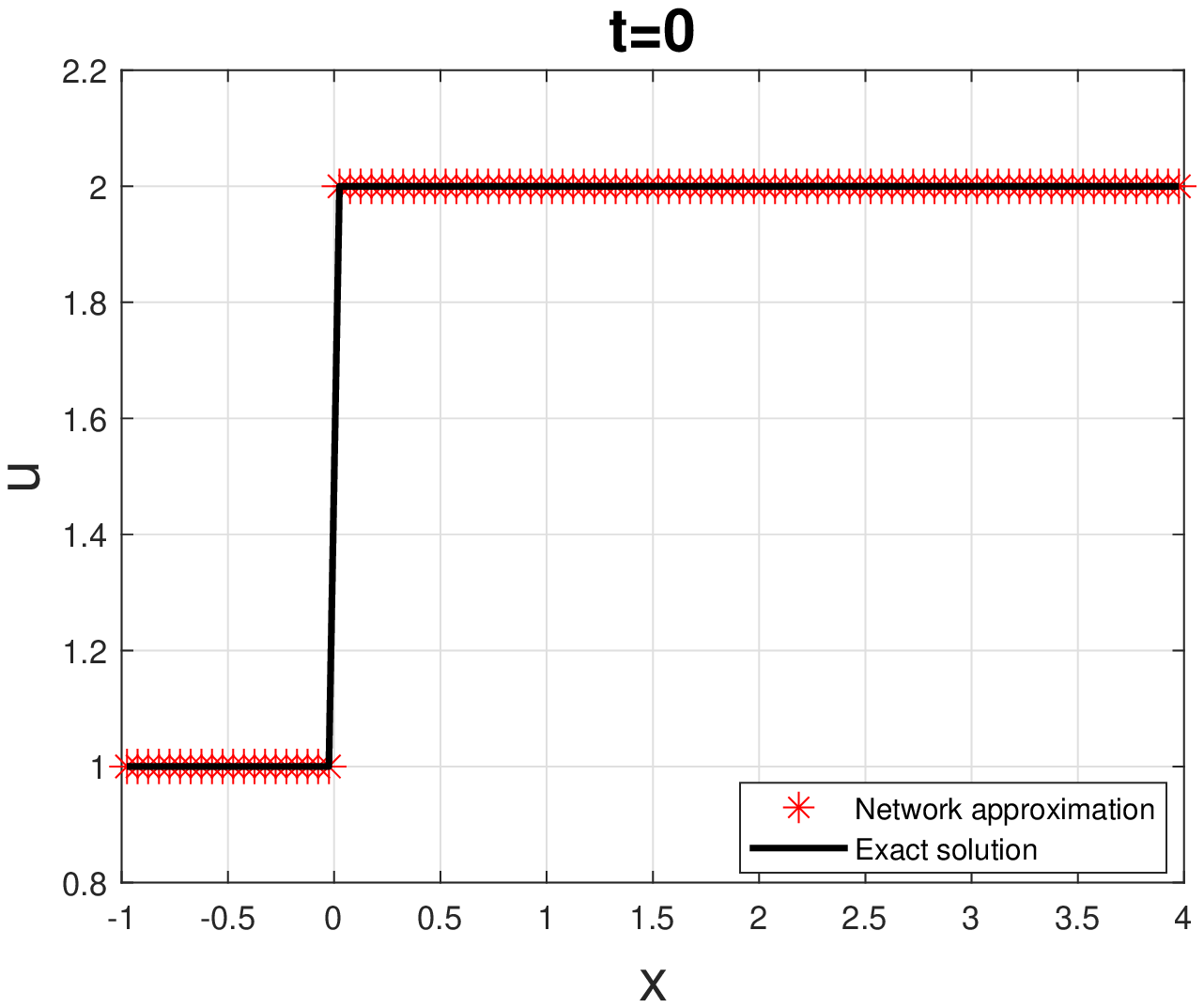}}\,
\subfigure[]{\includegraphics[width=0.24\linewidth]{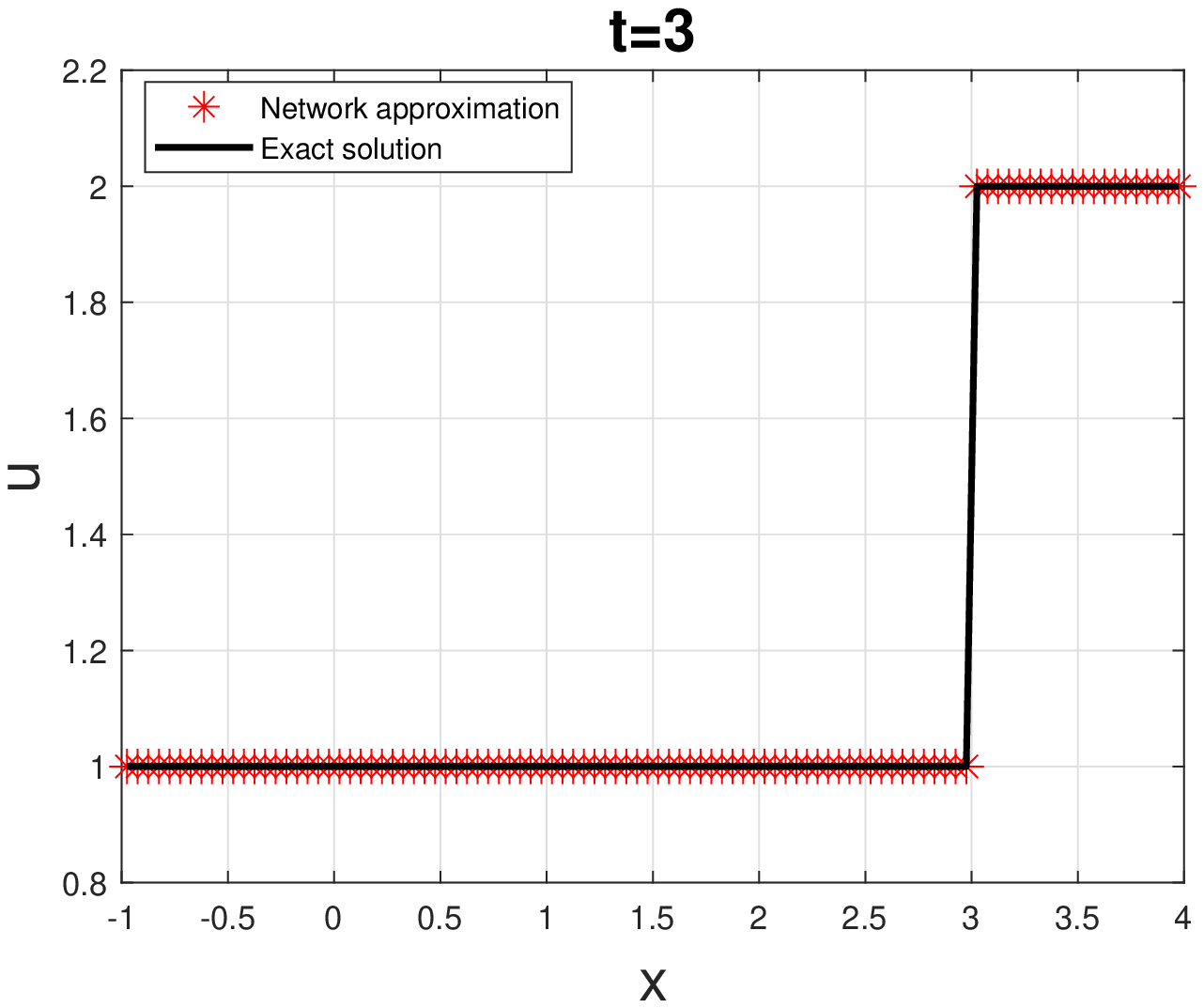}}
\caption{\small (a) \& (b) Inviscid Burgers' equation $u_t+(\frac{u^2}{2})_x=0$; (c) \& (d) Linear advection equation $u_t+u_x=0$}
\label{intro:fig-example}
\end{figure}

In a word, through intensive training we force the neural network to learn the solution average evolution between two neighboring time steps. The functionality of the parameter set $\Theta^{\star}$ is similar to the scheme definition of upwind finite volume method (\ref{intro:upwind-forward}), which are the coefficients of cell-average based neural network scheme (\ref{intro:CANN-scheme}). For nonlinear convection diffusion equation (\ref{intro:conv-diff}), same framework is applied that is summarized below
$$
\bar{u}_j(t_{n+1}) =\bar{u}_j(t_n)-\frac{1}{\Delta x}\int^{t_{n+1}}_{t_n}\int_{I_j}\left\{ f(u)_x-\mu u_{xx}\right\}~ dxdt \quad  \approx \quad
 \bar{v}_j^{n+1} = \bar{v}_j^{n}+\mathcal{N}(\vv_{j}^{in}; \Theta^{\star}).
$$
Notice that we do not discretize or approximate differential terms involving spatial variable. In stead, having the neural network handle all spatial variable related differentiation and integration approximation is the major idea of our method.
In Figure \ref{intro:fig-example} we list two simulations with our CANN method. Shock and contact discontinuity are both sharply captured.

Below we summarize the results and features of cell-average based neural network method. Some are quite outstanding and are not common to classical numerical methods.
\begin{itemize}
\item Adapt any time step size, i.e. $\Delta t=8\Delta x$, even being an explicit scheme
\item For Heat equation, errors are independent of time step size $\Delta t$.
 \item Introduce almost zero artificial numerical diffusion for contact discontinuity evolution
\end{itemize}
Due to some mysterious reason, the neural network method is able to catch up solution information around the next time level $t_{n+1}$ thus allows large time step size evolution. 
%Of course, it is fair to compare our neural network method to first order upwind scheme. The simplest network (i.e. 3 neurons with one hidden layer) have more than 10 parameters (or coefficients) than the regular finite volume or finite difference schemes, which can be tuned to maximize the performance. 
Different to classical numerical methods for which we design a scheme first, here the network itself finds the {\bf\emph{best scheme}} for us.

The organization of the article is the following. In section \S \ref{S:2.1}, we introduce the definition of cell-average based neural network method and highlight its connection to finite volume schemes. In section \S \ref{S:2.2}, we emphasize the training process and summarize the learning data difference between linear and nonlinear PDEs. In Section \S \ref{S:2.3}, we list the major results of cell-average based neural network method. Sequence of numerical examples are presented in section \S \ref{S:3}. Final conclusion remarks are given in section \S \ref{S:4}.

%%%%%%%%%%%%%%%%%%%%%%%%%%%%%%%%%%%%%%%%%%%%%%%%%%%%%%%%%%%%%%

\section{Neural network solver} \label{S:2}
\subsection{Problem setup, motivation and cell-averaged neural network method}
\label{S:2.1}
\vspace{.1in}

We consider to develop finite volume or cell-average based neural network (CANN) method solving partial differential equations (PDEs)
\begin{equation}\label{eq:PDE}
    u_t = \mathcal{L}(u),~~~(x, t) \in (a,b) \times \Re^{+}.
\end{equation}
%\begin{align}
%    &u_t = \mathcal{L}(u),~~~(x, t) \in \Omega \times \Re^{+}, \nonumber\\
%   &\mathcal{B}(u)=g, ~~~(x, t) \in \partial\Omega \times %Re^{+}, \nonumber\\
% &u(x,0)=u_0(x).
%\end{align}
Here $t$ and $x$ denote the time and spatial variables and $(a, b)$ is the spatial domain. Differentiation operator $\mathcal{L}$ is introduced to represent a generic first order hyperbolic or second order parabolic differentiation operator. For example, we have $\mathcal{L}(u)=(\frac{u^2}{2})_x$ for the inviscid Burgers equation and $\mathcal{L}(u)={u}_{xx}$ for the Heat equation. 

Our neural network method is mesh dependent and motivated by finite volume method. Once well trained. the cell-average based neural network method can be applied solving PDEs (\ref{eq:PDE}) as a regular finte volume scheme. Having a uniform partition of $(a, b)$ into $J$ cells and $\Delta x=\frac{b-a}{J}$ is adapted as the cell size. Denoting $x_{1/2}=a$, $x_{J+1/2}=b$, we have $[a, b]=\bigcup^J_{j=1} I_j$ with $I_j=[x_{j-1/2},x_{j+1/2}]$ as one computational cell. Furthermore, we have partition in time and adapt $\Delta t$ for the time step size and we have $t_n=n\times\Delta t$ with $t_0=0$.    

Now integrate the partial differential equation (\ref{eq:PDE}) over the  computational cell $I_j$ and time interval $[t_n, t_{n+1}]$, we have
\begin{equation}\label{eq:PDE-integration}
\int^{t_{n+1}}_{t_n}\int_{I_j} ~ u_t ~dxdt=\int^{t_{n+1}}_{t_n}\int_{I_j} ~ \mathcal{L}(u)~ dxdt.
\end{equation}
With the definition of cell average 
$\bar{u}_j(t)=\frac{1}{\Delta x}\int_{I_j} u(x,t) ~dx$, 
equation (\ref{eq:PDE-integration}) can be integrated out as
\begin{equation}\label{eq:PDE-integral-format}
\bar{u}_j(t_{n+1}) -\bar{u}_j(t_{n})=\frac{1}{\Delta x}\int^{t_{n+1}}_{t_n}\int_{I_j} ~ \mathcal{L}(u)~ dxdt.
\end{equation}
The equation above is nothing but the integral format or a weak formulation of (\ref{eq:PDE}). It is equivalent to the original partial differential equation (\ref{eq:PDE}) with extra differentiation condition applied. The integral format (\ref{eq:PDE-integral-format}) is the starting point at where we design our neural network method. The idea of cell-average based neural network method is to apply a simple fully connected network $\mathcal{N}(\cdot; \Theta)$ to approximate the right hand side of (\ref{eq:PDE-integral-format})
\begin{equation}\label{nn:network-goal}
    \mathcal{N}(\cdot; \Theta)\approx \frac{1}{\Delta x}\int^{t_{n+1}}_{t_n}\int_{I_j} ~ \mathcal{L}(u)~ dxdt,
\end{equation}
where $\Theta$ denotes the network parameter set of all weight matrices and biases. 

\begin{figure}[h]\label{nn:CANN-figure}
\centering
\includegraphics[width=0.4\linewidth]{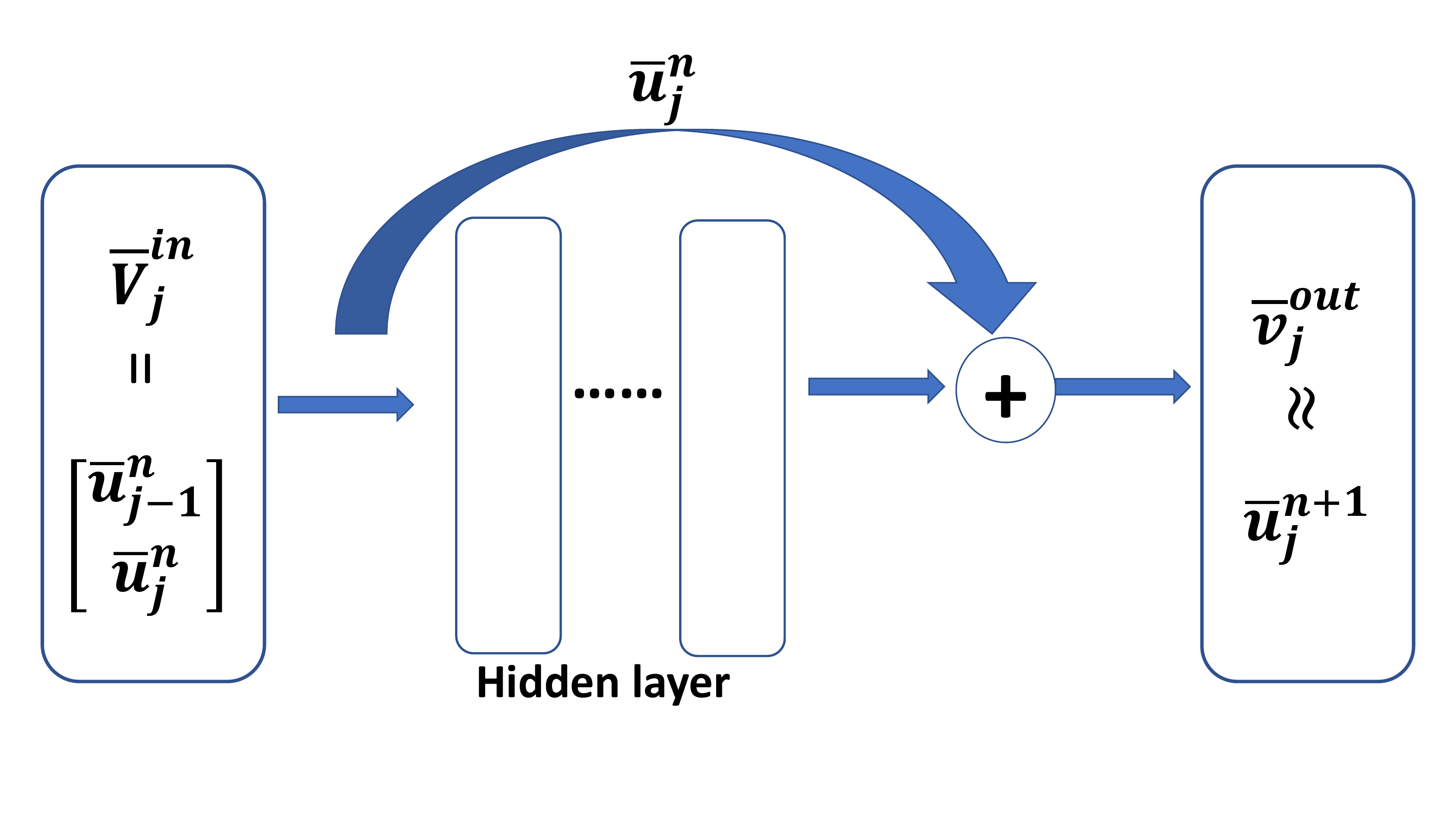}
\caption{Illustration of cell-average based neural network method}
\end{figure}

Given the solution averages $\left\{\bar{u}^n_j\right\}$ at time level $t_n$, we apply following neural network to approximate the solution average $\bar{u}^{n+1}_j$ at next time level $t_{n+1}$
\begin{equation}\label{nn:neural-network}
    \bar{v}_j^{out} = \bar{v}_j^{in}+\mathcal{N}(\vv_{j}^{in}; \Theta).
\end{equation}
With $\bar{v}_j^{in}=\bar{u}^{n}_j$ and comparing (\ref{nn:neural-network}) and the integral format (\ref{eq:PDE-integral-format}) of the PDEs, we have
$$
\bar{v}_j^{out}\approx \bar{u}^{n+1}_j.
$$
Vector $\vv_{j}^{in}$ as the input vector of the network $\mathcal{N}(\vv_{j}^{in}; \Theta)$ is an important component that should be carefully chosen, see Figure \ref{nn:CANN-figure}. Its general format is given as 
\begin{equation}\label{nn:input-vector-stencil}
    \vv_{j}^{in}=\Big[\bar{u}^n_{j-p}, \cdots, \bar{u}^n_{j-1}, \bar{u}^n_j, \bar{u}^n_{j+1}, \cdots,\bar{u}^n_{j+q} \Big]^T,
\end{equation}
%where we include $p$ cell averages to the left of $\bar{u}^n_j$ and $q$ cell averages to the right of $\bar{u}^n_j$ in the definition of the input vector. 
where we include the left $p$ cell averages and right $q$ cell averages of $\bar{u}^n_j$ in the input vector. 
The suitable stencil or the $p$ and $q$ values in (\ref{nn:input-vector-stencil}) determine the effectiveness of the neural network method approximating the solution average $\bar{u}_j^{n+1}$ at the next time level. 

Let's use advection equation $u_t+u_x=0$ as an example to demonstrate the relation between finite volume method and the cell-averaged neural network method. Forwar Euler upwind finite volume scheme can be rewritten in terms of vector multiplication format as
$$\bar{u}^{n+1}_j -\bar{u}^{n}_j=\frac{\Delta t}{\Delta x}\left(\bar{u}^{n}_j-\bar{u}^{n}_{j-1}\right)=
\begin{bmatrix} \frac{\Delta t}{\Delta x}& -\frac{\Delta t}{\Delta x}
\end{bmatrix}
\begin{bmatrix}
\bar{u}^{n}_j\\\bar{u}^{n}_{j-1}
\end{bmatrix}.
$$
Now consider the cell-average based neural network method for $u_t+u_x=0$. We adapt upwind mechanism and choose input vector as $\vv_{j}^{in}=[\bar{u}^{n}_j,\bar{u}^{n}_{j-1}]^T$, see illustration in Figure \ref{nn:CANN-figure}. The goal of the neural network is to obtain optimal parameter set $\Theta$ through minimizing the error between network output $\bar{v}_j^{out}$ and the given target $\bar{u}^{n+1}_j$. The parameter set $\Theta$ is similar to the coefficients of a finite volume scheme, for example the $\frac{\Delta t}{\Delta x}$ before $\bar{u}^{n}_j$ and $-\frac{\Delta t}{\Delta x}$ before $\bar{u}^{n}_{j-1}$ of the upwind scheme. We highlight the network parameter set $\Theta$ interact with input vector $\vv_{j}^{in}$ in a nonlinear fashion, even for a linear differential equation. Thus neural network method behaves differently to finite volume schemes. Once well trained, the optimal parameter set $\Theta^{*}$ uniquely determines the {\bf\emph{scheme definition}} of a neural network method. We then apply this neural network solver at any cell location (index $j$) and at any time level $t_n$ to evolve the solution average to next time level $t_{n+1}$ as an explicit scheme.

Motivated by the well established finite volume schemes for hyperbolic PDEs, we adapt upwind or method of characteristic mechanism to choose suitable network input vector $\vv_{j}^{in}$ for hyperbolic PDEs. And we adapt a symmetric mechanism similar to central scheme for Heat equation to choose suitable network input vector for parabolic PDEs. The guideline of picking a suitable stencil ($p$ and $q$ in (\ref{nn:input-vector-stencil})) for the network input vector $\vv_{j}^{in}$ is summarized below. 

\vspace{.1in}
\noindent{\bf Guideline on the choice of network input vector $\vv_{j}^{in}$:}
\begin{itemize}
    \item Hyperbolic PDEs: include the characteristic or domain of dependence 
    \item Parabolic PDEs: a symmetric stencil mechanism
\end{itemize}
In this paper, we consider a standard fully connected neural network with $M$ ($M\geq 3$) layers. The input and output vectors of the network are the first and last layer. Among the total $M$ layers, the interior $(M-2)$ are the hidden layers. Thus the minimum structure of the neural network involves one hidden layer with $M=3$. We have $n_i$ ($i =1, \cdots, M$) denote the number of neurons in each layer. The first layer is the input vector with $n_1=p+q+1$ as its dimension and the last layer is the output vector with $n_M=1$ as its dimension. The abstract goal of machine learning is to find a function $\mathcal{N}: R^{p+q+1} \to R^1$ such that $\mathcal{N}(\cdot ; \Theta)$ accurately approximates $\frac{1}{\Delta x}\int^{t_{n+1}}_{t_n}\int_{I_j} ~ \mathcal{L}(u)~ dxdt$, the right hand side of (\ref{eq:PDE-integral-format}). The optimal parameter set $\Theta$ of the network $\mathcal{N}(\cdot ; \Theta)$ will be obtained by training the network intensively over the given data set.

Learning data are collected in the form of pairs. Each pair refers to the solution averages at two neighboring time level. The training data set is denoted as
\begin{equation}\label{nn:learning-data-set}
S =\left\{\left(\bar{u}_{j}^n,\bar{u}_{j}^{n+1}\right),\,\,j=1, \cdots, J\right\}_{n=0}^{m},
\end{equation}
which are solution averages obtained from highly accurate numerical method for the PDE. We highlight the training data pairs are solution averages collected over the spatial domain and from time levels $t_0$ to $t_{m+1}$. 

The fully neural network $\mathcal{N}(\cdot; \Theta)$ consists of $M$ layers. Each two consecutive layers is connected with an affine linear transformation and a point-wise nonlinear activation function.
The function or the mapping $\mathcal{N}(\cdot; \Theta)$ is a composition of following operators,
\begin{equation}
    \mathcal{N}(\cdot; \Theta)=(\sigma_M\circ \vW_{M-1})\circ \cdots \circ (\sigma_2\circ \vW_1). \label{neuralNetwork-structure}
\end{equation}
where $\circ$ stands for operator composition. We have $\bf{W}_i$ denoting the linear transformation operator or the weight matrix connecting the neurons from $i$-th layer to $(i+1)$-th layer. The parameter set is further augmented with the biases vectors. We have $\sigma_i: R\to R$ as the activation function ($i\geq 2$), which is applied to each neuron of the $i$-th layer in a component-wise fashion. In this paper we apply $tanh(x)$ function as the activation function. Specifically $\sigma_i=tanh(x)$ is applied between all layers, except to the output layer for which we have $\sigma_M(x) =x$.

\begin{definition}\label{Def:1}
A cell-average based neural network method is uniquely determined by the following four components: (1) the choice of spatial mesh size $\Delta x$; (2) the choice of time step size $\Delta t$; (3) the choice of network input vector $\vv_j^{in}$ of (\ref{nn:input-vector-stencil}); and (4) the number of hidden layers and neurons per layer the corresponding structure of the neural network. 
\end{definition}

We further highlight that our cell-average based neural network method, once well trained, will be implemented as a regular explicit finite volume scheme. Notice the cell-average based neural network method is designed to approximate the weak or integral format of the partial differential equations (\ref{eq:PDE-integral-format}), not the original partial differential equations (\ref{eq:PDE}) given in differentiation format.

%%%%%%%%%%%%%%%%%%%%%%%%%%%%%%%%%%%%%%%%%%%%%%%%%%%%%%%%%%%%%%%%%%
\subsection{Training Process}
\label{S:2.2}
\vspace{.1in}

In this section, we discuss how to train the network to obtain optimal parameter set $\Theta^{*}$ such that the neural network (\ref{nn:neural-network}) can accurately approximate the solution average evolution $\uAve^n_{j}\to \uAve^{n+1}_j$. This is achieved by applying  $\bar{v}^{in}_j=\uAve_j^n$ in (\ref{nn:neural-network}) to obtain  network output $\bar{v}_j^{out}$, comparing with the target $\uAve_j^{n+1}$, and then looping among the data set $S$ of (\ref{nn:learning-data-set}) to minimize the error or the squared loss function
\begin{equation}\label{nn:loss-function}
    L_{j,t_n}(\Theta)=(\bar{v}_j^{out}-\uAve_j^{n+1})^2, 
\end{equation}
for all $j =1,\cdots,J$ and for all $n=0,\cdots, m$. This choice of loss function defined over one single data pair is referred as the stochastic or approximate gradient descent method. Notice for $j=1$ or $j=J$ or those close to boundary cells, the stencil or network input vector $\vv_j^{in}$ of (\ref{nn:input-vector-stencil}) requires averages values of $p$ cells to the left and $q$ cells to the right of the current cell. In this paper we only consider Dirichlet or periodic boundary conditions. Thus we simply copy solution averages from inside the domain implementing periodic boundary condition. And we assign exact solution values to those out of domain ghost cells solving Burgers equation Riemann problems.

Recall we have $S$ of (\ref{nn:learning-data-set}) denoting the training data set, which include solution averages spread over the spatial domain with index $j$ and up to time levels $t_{m+1}$. Numerical tests show that for linear partial differential equations,  training data of $(t_0, t_1)$ with $m=0$ in $S$ is sufficient for obtaining optimal parameter set $\Theta$. For nonlinear partial differential equations, i.e. the Burgers equation, multiple time levels ($m>0$ in $S$) data are necessary to guarantee the neural network effectively learning the solution evolution mechanism, for example capturing the right shock speed.

\begin{figure}[h]
\centering
\includegraphics[width=0.55\linewidth]{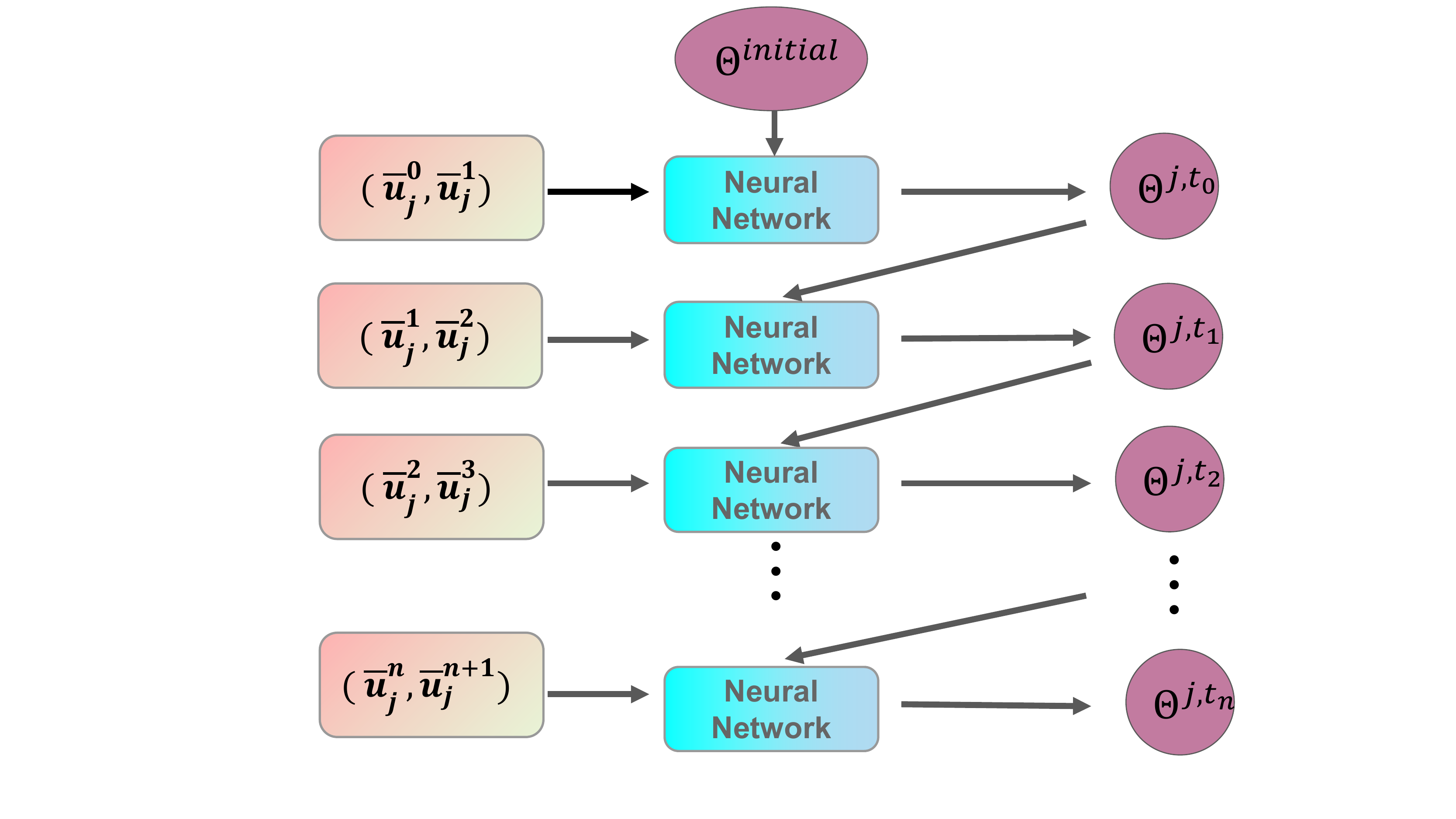}
\caption{Network parameter one epoch training for cell-average based neural network (CANN) method}\label{nn:training-order}
\end{figure}

As discussed in section \ref{S:2.1}, the network parameter set $\Theta$ behave more or less as the coefficients of a finite volume scheme. Parameter set $\Theta$ is independent of at what cell location and at what time level the neural network method is applied. Thus we should minimize the loss function over all spatial index $j$ and at all time levels in $S$ to obtain the optimal parameter set $\Theta^{*}$ or the {\emph{coefficients}} of the neural network method. 

One epoch of iteration is defined as the following. The weights and biases are first randomly generated from normal distribution around zero. The parameter set $\Theta$ are updated sequentially through space and time with all data pairs applied {\emph{once}} in the training set $S$. In a word, we start with first time level pair $(t_0, t_{1})$ and go through all spatial cells and iteratively update the parameter set $\Theta^{j,t_0}$ through the gradient descent direction of (\ref{nn:loss-function}). Then we move on to the next time level pairs $(t_n, t_{n+1})$ ($n\geq 1$) and go through all spatial cells again updating the parameter set $\Theta^{j,t_n}$. Below is an illustration of one epoch definition, for which we use solution averages pairs in the training set $S$ once iterating the network parameter, see Figure \ref{nn:training-order}.
\begin{equation*}\label{NN:one-epoch}
    \Theta^{j,t_n}_{old}\longrightarrow \Theta^{j,t_n}_{new},\quad j=1,\cdots, J,\quad n=0,\cdots, m.
\end{equation*}
The superscripts $j$ and $t_n$ notations are added to specify the iteration or update is processed at what cell and at what time level. They all refer to the update of the same optimized parameter set $\Theta^{*}$. Learning rate $\alpha=0.01$ or $\alpha=0.001$ is applied, when updating through stochastic gradient descent method. 

In this paper, we have integer $K$ introduced denoting the total number of epochs iterations involved in training. Instead of iterating through one epoch defined above, at each time level we run total $K$ iterations to update the parameter set $\Theta$ with one iteration defined as the one round of updating over spatial domain (index $j$). We then continue the iteration of $\Theta$ to the next time level, which can be summarized as
\begin{equation}\label{nn:theta-iteration}
    \Theta^{j,t_n,i}\longrightarrow \Theta^{j,t_n,i+1},\quad j=1,\cdots, J,\quad i=1,\cdots, K, \quad n=0,\cdots, m.
\end{equation}
In the end, the parameter set $\Theta$ have been iteratively updated with total $J\times K\times (m+1)$ times. 
Specifically for the last time level pair $(t_m,t_{m+1})$, we record the squared $L_2$ error defined below corresponding to iteration
\begin{equation}\label{nn:L2-error-square}
   L^2_2(t_{m+1})=\sum^J_{j=1}  \Delta x L_{j,t_m}(\Theta)=\sum^J_{j=1}(\bar{v}_j^{out}-\uAve_j^{m+1})^2\Delta x.
\end{equation}
We output the squared $L_2$ error of \eqref{nn:L2-error-square} corresponding to iteration index $i=1,\cdots,K$ to demonstrate the effectiveness of cell-average based neural network method. %capturing the solution evolution mechanism. 
Now we conclude the section with comment on the minimum size of training data set $S$ of (\ref{nn:learning-data-set}), which are purely lab results observed from numerical tests.

\begin{rem}{\label{remark:1}}
For linear partial differential equations, one time level solution averages in the training set $S$ corresponding to $(t_0,t_1)$ or $m=0$ in (\ref{nn:learning-data-set}), is sufficient for obtaining an effective neural network. For nonlinear partial differential equations, it is necessary to include multiple time levels ($m>0$) of solution averages in the training set $S$ to have the neural network learn the evolution mechanism successfully.
\end{rem}

For example in the numerical section when we solve inviscid Burgers' equation with smooth $\sin(x)$ initial in Example {\ref{ex5}}, we have solution averages data pairs up to $t=2$ included in the training set $S$. This is because a shock starts to develop at $t=1$. We need to have more time levels solution averages included in the training set to make sure the neural network solver be able to learn the shock capturing mechanism.  

%%%%%%%%%%%%%%%%%%%%%%%%%%%%%%%%%%%%%%%%%%%%%%%%%%%%%%%%%%%%%%%%%
\subsection{Implementation and summary of cell-average based neural network method}
\label{S:2.3}
\vspace{.1in}

With the optimal weights and biases $\Theta^{*}$ obtained and the neural network $\mathcal{N}(\vv_{j}^{in}; \Theta^{*})$ well defined and available, the cell-average based neural network method can be implemented as a regular {\bf \emph{explicit finite volume scheme}} as below
\begin{equation}\label{nn:neural-scheme}
    \bar{v}_j^{n+1} = \bar{v}_j^{n}+\mathcal{N}(\vv_{j}^{n}; \Theta^{*}),\quad \forall j=1,\cdots,J, \quad \forall n=0,1,2,\cdots.
\end{equation}
Again, with the spatial and time step sizes of $\Delta x$ and $\Delta t$ and the previously chosen network input vector of  $\vv_{j}^{n}=\Big[\bar{v}^n_{j-p}, \cdots, \bar{v}^n_{j-1}, \bar{v}^n_j, \bar{v}^n_{j+1}, \cdots,\bar{v}^n_{j+q} \Big]^T$, together with the network optimal parameter set $\Theta^{*}$, we have a complete definition of a neural network method.

For cell-average based neural network method (\ref{nn:neural-scheme}), we assign ghost cell values and apply boundary conditions as a regular finite volume method. Again, we only consider Dirichlet or periodic boundary conditions in this paper. A generic stencil include $p$ cells to the left and $q$ cells to the right to evolve current cell to the next time level. Similar to the training process discussed in section \ref{S:2.2}, we simply copy solution averages from inside the domain for periodic boundary conditions. And we assign exact solution values to those out of domain ghost cells for Burgers' equation Riemann problems, i.e. for Dirichlet boundary conditions.

One amazing result is that cell-average based neural network method can be {\bf\emph{relieved}} from the CFL restriction on time step size, especially for parabolic PDEs. With explicit discretization in time, classical numerical method requires time step size to be as small as $\Delta t\approx (\Delta x)^2$, which is very expensive for multi-dimensional problems. It turns out neural network method can adapt to any time step size $\Delta t$ that is independent of spatial size $\Delta x$, and still a stable method can be obtained. Once well trained, the neural network method can efficiently and accurately evolve the solution forward in time as an explicit method.

Recall our neural network method is based on the following integral format of the partial differential equation $u_t+\mathcal{L}(u)=0$
\begin{equation*}
\bar{u}_j(t_{n+1}) -\bar{u}_j(t_n)=\frac{1}{\Delta x}\int^{t_{n+1}}_{t_n}\int_{I_j} ~ \mathcal{L}(u)~ dxdt,
\end{equation*}
with $\mathcal{L}$ as the differentiation operator on the spatial variable. It is equivalent to the PDEs given in differentiation format. Due to some mysterious reason, the neural network method is able to catch up solution information around time level $t_{n+1}$ thus allows large time step size evolution as an implicit method. Below we summarize the major result of our CANN method. 
\begin{lemma}
Even being an explicit scheme, cell-average based neural network method (\ref{nn:neural-scheme}) can adapt to any time step size $\Delta t$ for  hyperbolic and parabolic PDEs (\ref{eq:PDE}). With $\Delta t=\Delta x$, first order of accuracy is obtained with neural network method (\ref{nn:neural-scheme}) for linear advection and linear convection diffusion equations. For Heat equation, neural network method errors depend only on spatial mesh size $\Delta x$ and are independent of time step size $\Delta t$.
\end{lemma}

Besides accuracy and error behavior studies, we also carry our a series of numerical tests with cell-average based neural network method. Now we list the advantages of our CANN method that are not common for classical numerical methods.

\begin{itemize}
  \item Allow large time step size evolution, i.e. $\Delta t=8\Delta x$ for linear hyperbolic PDEs.
  \item For Heat equation, errors are independent of time step size. Similar errors are obtained for $\Delta t=\Delta x$, $\Delta t= 2\Delta x$ and $\Delta t=4\Delta x$ when same spatial mesh $\Delta x$ is adapted.
  \item Introduce almost zero artificial numerical diffusion for contact discontinuity propagation.
  \item Introduce little numerical dissipation and dispersive errors after long time run.
\end{itemize}

\begin{rem}
Once one cell-average based neural network is well trained and available, it can be applied to solve same PDE associated with different initials and over different domains.
\end{rem}

\begin{rem}
It remains unknown how to choose a suitable neural network architecture in terms of number of hidden layers and neurons per layer. Numerical tests show one or two hidden layers with a few neurons work well. We further mention that some network structure may lead to extremely small errors and can be even regarded as a perfect solver, for which regular order of convergence over refined mesh error analysis does not hold anymore. 
\end{rem}

%%%%%%%%%%%%%%%%%%%%%%%%%%%%%%%%%%%%%%%%%%%%%%%%%%%%%%%%%%%%%%%%%%%%%%%%
\section{Numerical Example}
\label{S:3}
In this section, we carry out a series of numerical tests to check out the accuracy and capability of cell-average based neural network methods. We start with linear advection equation, Heat equation and convection diffusion equations to investigate if or not order of convergence can be observed. Then we move on to nonlinear hyperbolic conservation law to test the capability of neural network method capturing shock and rarefaction waves propagation.

Over the section we adapt $T$ as a generic final time at where we compute the errors and orders. As mentioned previously, time step size $\Delta t$ is chosen before and is a part of the definition of a neural network method (\ref{nn:neural-scheme}). Thus we always have $T$ as an integer multiple times of $\Delta t$. Below we list the $L_2$ and $L_{\infty}$ errors formula as used in finite volume methods.

\begin{equation}\label{nn:L2-norm-error}
  Error_{L_2}(T)=\sqrt{\sum^J_{j=1}(\bar{v}_j(T)-\uAve_j(T))^2\Delta x}
\end{equation}
\begin{equation}\label{nn:Linfty-norm-error}
  Error_{L_{\infty}}(T)=\max^J_{j=1}|\bar{v}_j(T)-\uAve_j(T)|
\end{equation}
Again, we have $\bar{v}_j(T)$ denote our CANN method (\ref{nn:neural-scheme}) solution on cell $j$ and at final time $T$. We have $\bar{u}_j(T)$ denote either the exact solution average or the reference solution average on cell $j$ and at time $T$ that is obtained from a highly accurate numerical method. 

As marked down in Definition \ref{Def:1}, four components of $\Delta x$, $\Delta t$, network input vector $\vv_j^{in}$ and the structure of neural network (number of hidden layers and neurons per layer) together identify one neural network solver (\ref{nn:neural-scheme}). Most of the time we follow the principle listed in section {\ref{S:2.1}} to choose suitable stencil width or the network input vector of $\vv_j^{in}$. We highlight that for linear PDEs one time level $(t_0, t_1)$ solution average data pair is applied training the network and for nonlinear PDEs multiple time levels solution average data pair are needed to obtain an effective neural network solver, see Remark {\ref{remark:1}}. We also mention that squared $L_2$ error of (\ref{nn:L2-error-square}) around $10^{-8}$ or smaller is used as stop condition in training.

\subsection{Linear advection equation}\label{S3-1}
In this subsection, we focus on linear hyperbolic equation of
\begin{equation}\label{eq:transport}
u_t+u_x=0.
\end{equation}
Even the above linear equation is a simple model, quite a few problems can be tested to evaluate the capability of a numerical method. We consider three problems for advection (\ref{eq:transport}). One is about smooth function evolution and we check if order of convergence can be observed. Then we study the contact discontinuity propagation problem. It is not a trivial test, since numerical methods tend to either generate smeared out approximations or numerical oscillations around the discontinuity. For the third test, we have the neural network method simulate wave propagation after long time run. Neural network method is quite stable and produces little dispersive and dissipation errors after long time simulation.

\begin{example}\label{ex1}{\bf \emph {smooth wave propagation}}
\end{example}

In this example we solve (\ref{eq:transport}) with initial condition $u(x,0)=\sin(x)$. Spatial domain is taken as $D=[0, 2\pi]$ and the exact solution $u(x,t)=\sin(x-t)$ is a smooth wave propagating from left to right.
We consider a simple neural network with input vector of 
\begin{equation}\label{eq:upwind-input-vector}
    \vv_{j}^{in} = \Big[ \uAve^n_{j-1}, \uAve^n_j\Big]^T,
\end{equation}
that is similar to the upwind finite volume scheme. The network picked consists of 2 hidden layers with 5 neurons per layer. Total iterations of $K=5\times 10^5$ is applied to train the networks. Periodic boundary conditions are applied. We carry out two accuracy tests, for which we compute the $L_2$ and $L_{\infty}$ errors of (\ref{nn:L2-norm-error}) and (\ref{nn:Linfty-norm-error}) at final time $T=\pi$.

\vspace{.1in}
\noindent
{\bf Case I:}

For the first case study we mimic the accuracy check of a standard numerical method. We consider four spatial meshes of $\Delta x=\pi/10$, $\pi/20, \pi/40,$ and $ \pi/80$. For each cell size $\Delta x$, time step size $\Delta t=\Delta x$ is taken correspondingly. The four well trained neural network solvers are very similar to each other except the mesh size. All four network solvers are used to solve the smooth wave propagation to final time $T=\pi$ and $L_2$ and $L_{\infty}$ errors are computed and listed in Table \ref{test1-table1}. 
%Roughly first order of accuracy is obtained. We mention that the errors tend to cease decreasing on refined mesh ($\Delta x=\pi/160$).

\begin{table}[!htb]
\centering
\begin{tabular}[c]{l l l l l l}
\hline
$\Delta x$ & $L_2$ & order & $L_{\infty}$ & order\\
\hline
$\pi/10$  & $1.8756e^{-2}$ & & $1.0237e^{-2}$&\\
$\pi/20$  & $8.0830e^{-3}$ &1.21 & $4.7403e^{-3}$&1.11\\
$\pi/40$  & $1.5547e^{-3}$ &2.39 & $9.7037e^{-4}$&2.29\\
$\pi/80$  & $6.3500e^{-4}$ &1.28 & $3.9838e^{-4}$ &1.28\\
\hline
\end{tabular}
\caption{Errors and orders of neural network methods for Example \ref{ex1} (Case I), $\Delta t=\Delta x$}\label{test1-table1}
\end{table}

\vspace{.1in}
\noindent
{\bf Case II:}

For this case we keep the spatial mesh size $\Delta x=\pi/40$ fixed, but consider three different time step sizes of $\Delta t=2\Delta x$, $\Delta t=5\Delta x$ and $\Delta t=8\Delta x$. The three neural network solvers apply same upwind like network input vector of (\ref{eq:upwind-input-vector}). The only difference is the time step size. We list the $L_2$ and $L_{\infty}$ errors computed at final time $T=\pi$ in Table \ref{test1-table2}. Three network solvers give similar errors and it is not clear if the network errors relate to $\Delta t$. For all three time step sizes, i.e. $\Delta t=8\Delta x$, the principle of having domain of dependence included is not followed. We simply apply (\ref{eq:upwind-input-vector}) as the input vector that is like the upwind scheme. These settings conflict with method of characteristic. But all neural network solvers work well and give errors similar to those in case I.  
%Accuracy of cell-average based neural network method seems to be independent of the time step size.

\begin{table}[!htb]
\centering
\begin{tabular}[c]{l l l l}
\hline
$\Delta t$ & $L_2$ & $L_{\infty}$ \\
\hline
2$\Delta x$ & $7.0431e^{-3}$ & $4.1544e^{-3}$\\
5$\Delta x$ & $9.2344e^{-3}$ & $5.3628e^{-3}$\\
8$\Delta x$ & $6.9895e^{-3}$ & $3.1796e^{-3}$\\
\hline
\end{tabular}
\caption{Errors and orders of neural network methods for Example \ref{ex1} (Case II), $\Delta  x=\pi/40$ fixed with varying $\Delta  t$}\label{test1-table2}
\end{table}

\begin{figure}[!htb]
\centering
\includegraphics[width=0.3\linewidth]{fig/example2/u_0.eps}\quad\quad
\includegraphics[width=0.3\linewidth]{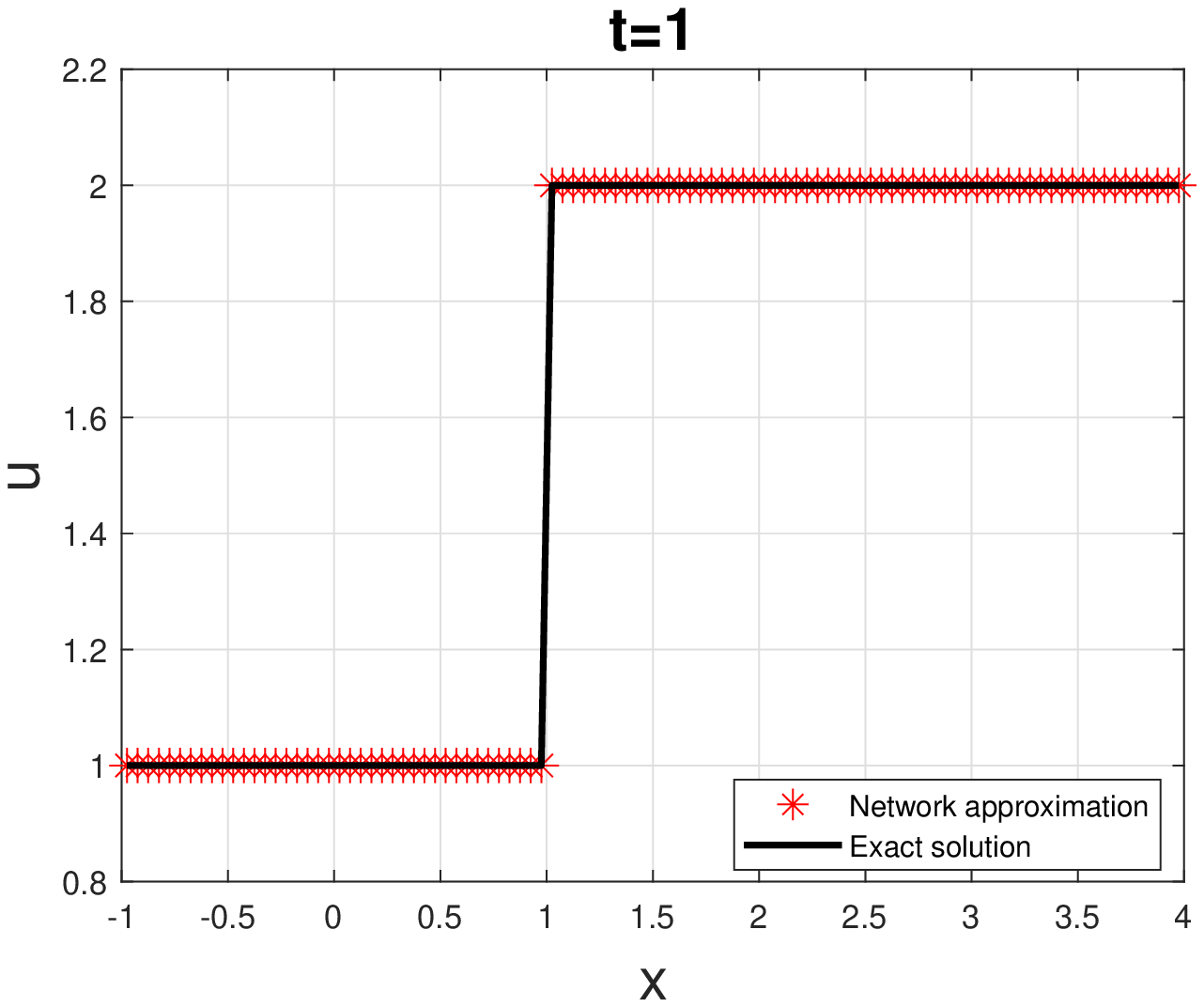}\\
\includegraphics[width=0.3\linewidth]{fig/example2/u_3.eps}\quad\quad
\includegraphics[width=0.3\linewidth]{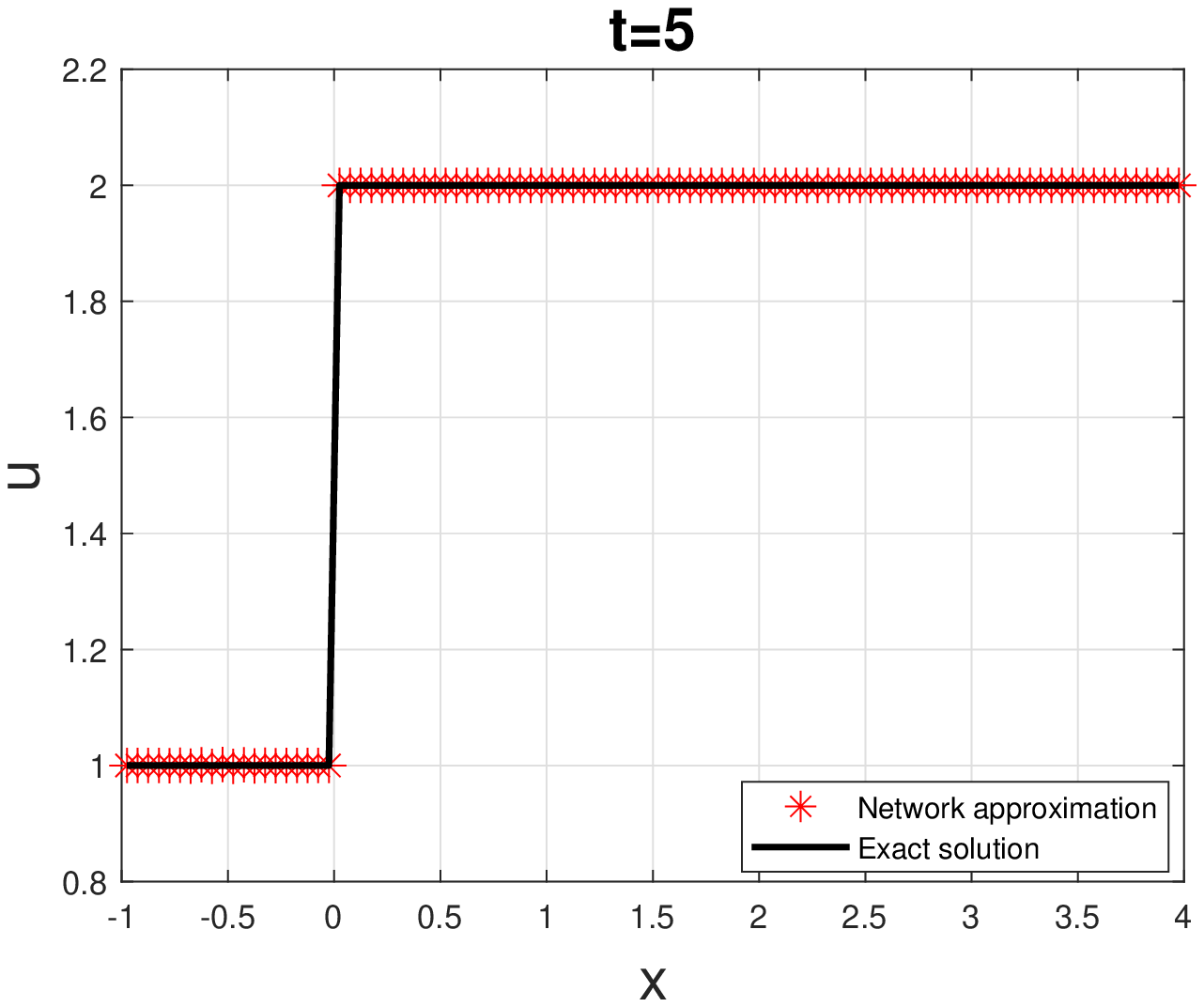}
\caption{contact discontinuity evolution (Example \ref{ex2}) with neural network method}
\label{fig:contact-discontinuity}
\end{figure}

\begin{example}\label{ex2}{\bf \emph {contact discontinuity}}
\end{example}
In this example, we solve advection equation (\ref{eq:transport}) with initial condition 
$$u(x,0)=\left\{
\begin{aligned}
1 ~~ x\le 0, \\
2 ~~ x> 0,
\end{aligned}
\right. $$
over domain $D=[-1, 4]$. Periodic boundary condition is applied. The contact discontinuity is initially located at $x=0$, which moves back into the domain after $t>4$. Same network input vector of (\ref{eq:upwind-input-vector}) that is similar to upwind scheme is considered. One hidden layer of 10 neurons is the chosen network structure and a total of $K=5\times 10^5$ iterations are applied training the network. We have $\Delta x= \Delta t=\frac{1}{20}$. Snapshots of our CANN method simulation are presented in Figure \ref{fig:contact-discontinuity}.
The contact discontinuity is sharply resolved, even after one period of evolution at $t=5$. No oscillation is observed and there is almost no artificial diffusion introduced with the neural network method. Numerical simulation with CANN method behaves better than many numerical methods that tend to generate smeared out simulations. 

\begin{example}\label{ex3}{\bf\emph {long time evolution}}
\end{example}
In this example we solve advection equation (\ref{eq:transport}) with initial condition 
$$
u(x,0)=\sin(x),
$$ over domain $D=[0, 2\pi]$. We choose $\Delta x= \frac{2\pi}{100}$ and $\Delta t=4\Delta x$. Long time simulation with CANN method is considered. To include the domain of dependence, following neural network input vector is taken  
\begin{equation}\label{eq:upwind-input-vector}
    \vv_{j}^{in} = \Big[ \uAve^n_{j-6},\uAve^n_{j-5},\uAve^n_{j-4},\uAve^n_{j-3},\uAve^n_{j-2},\uAve^n_{j-1}, \uAve^n_j\Big]^T.
\end{equation}
The chosen network is with 1 hidden layer and 10 neurons. A total of $K=5\times10^5$ iterations are applied. We compute $L_2$ and $L_{\infty}$ errors at different time locations, which are listed in Table \ref{test3-table}. After several periods, neural network method is still able to accurately capture the solution evolution and gives small dispersive and dissipation errors.
\begin{table}[!htb]
\centering
\begin{tabular}[c]{l l l l l}
\hline
& $t=4\pi/5$ & $t=2\pi$ & $t=4\pi$ & $t=8\pi$\\
\hline
$L_2$  & $9.1869e^{-14}$ & $1.0045e^{-12}$ & $1.8625e^{-10}$ & $7.1709e^{-6}$\\
\hline
$L_{\infty}$ & $1.2620e^{-13}$ & $1.0600e^{-12}$ & $1.7158e^{-10}$ & $5.3325e^{-6}$\\
\hline
\end{tabular}
\caption{$L_2$ and $L_{\infty}$ errors (Example \ref{ex3}) after long time simulation with neural network method}\label{test3-table}
\end{table}

% \begin{figure}[!htb]\label{test-3}
% \begin{minipage}[c]{0.40\linewidth}
% \centering
% \includegraphics[width=1\linewidth]{fig/test3_fig.eps}
% \caption{Trajectory approximations for different long time run}\label{test3-figure}
% \end{minipage}
% \begin{minipage}[c]{0.40\linewidth}
% \centering
% \begin{tabular}[c]{l l l l l}
% \hline
% & $t=4\pi/5$ & $t=7.5$ & $t=4\pi$ & $t=7.2\pi$\\
% \hline
% $L_2$  & $3.9258e^{-4}$ & $6.0389e^{-4}$ & $8.2124e^{-4}$ & $5.1561e^{-3}$\\
% \hline
% $L_{\infty}$ & $2.3403e^{-4}$ & $3.7759e^{-4}$ & $5.2533e^{-4}$ & $4.9877e^{-3}$\\
% \hline
% \end{tabular}
% \caption{$L_2$ error and $L_{\infty}$ error for different long time simulations}\label{test3-table}
% \end{minipage}
% \end{figure}

\subsection{Linear convection diffusion equation}
\label{S3-2}

In this section, we solve linear convection diffusion equation with cell-average based neural network method. It is further confirmed that we can be relieved from explicit scheme small time step size restriction. For Heat equation, numerical tests show the errors only relate to spatial mesh size $\Delta x$ and is independent of time step size $\Delta t$.

\begin{example}\label{ex4} {\bf \emph {Heat equation}}
\end{example}
In this example, we consider solving the Heat equation 
\begin{equation*}%\label{eq:Heat-equation}
u_t= u_{xx},
\end{equation*}
with initial $u(x,0)=\sin(\pi x)$ and over domain $D=[0, 1]$. Exact solution is $u(x,t)=e^{-\pi^2 t}\sin(\pi x)$. Different to hyperbolic PDEs, we have infinite speed of propagation for Heat equation. Here we follow a  similar to central scheme symmetric mechanism to pick up the network input vector as
\begin{equation}\label{eq:heat-input-vector}
    \vv_{j}^{in} = \Big[\uAve^n_{j-3}, \uAve^n_{j-2}, \uAve^n_{j-1}, \uAve^n_j, \uAve^n_{j+1}, \uAve^n_{j+2},\uAve^n_{j+3}\Big]^T.
\end{equation}
The chosen network is with 2 hidden layers and 15 neurons per layer. Iterations are run up to $K=10^5$ times to optimize the network parameter set. Periodic boundary conditions are applied. Final time $T=0.1$ is adapted to compute $L_2$ and $L_{\infty}$ errors of (\ref{nn:L2-norm-error}) and (\ref{nn:Linfty-norm-error}). We carry out three accuracy tests. For the first and third tests, we always choose $\Delta t=\Delta x$ and have $\Delta x$ refined to check out order of convergence. In the second test we have $\Delta x$ fixed and choose $\Delta t$ as a multiple times of $\Delta x$. Same network structure is applied, but different choices of $\Delta t$, $\Delta x$ or network input vector $\vv_j^{in}$ are considered in each test.

\vspace{.05in}
\noindent
{\bf Case I:}

For this case, we investigate the convergence order of neural network method. We consider four settings of $\Delta x= \frac{1}{40}, \frac{1}{80}, \frac{1}{160}, \frac{1}{320}$. For each cell size $\Delta x$, we choose $\Delta t=\Delta x$ correspondingly. The four network solvers are similar to each other except the mesh size. We solve the Heat equation with each neural network solver to $T=0.1$ and compute the $L_2$ and $L_{\infty}$ errors. In Table \ref{test4-table4} we list all errors and orders. Clean first order of convergence is observed with CANN method.

\begin{table}[!htb]
\centering
\begin{tabular}[c]{l l l l l l}
\hline
$\Delta x$ & $L_2$ & order & $L_{\infty}$ & order\\
\hline
$1/40$  & $8.6949e^{-3}$ & & $2.0873e^{-2}$&\\
$1/80$ & $4.5270e^{-3}$ & 0.94& $1.4104e^{-2}$ &0.57\\
$1/160$ & $2.4736e^{-3}$ &0.87 & $7.2650e^{-3}$ &0.96\\
$1/320$ & $1.2894e^{-3}$ &0.94 & $3.7860e^{-3}$ &0.94\\
\hline
\end{tabular}
\caption{Errors and orders of neural network methods for Heat equation (Case I), $\Delta t=\Delta x$}\label{test4-table4}
\end{table}

\vspace{.05in}
\noindent
{\bf Case II:}

For this case we fix the spatial mesh size $\Delta x=1/160$ but apply three time step sizes of $\Delta t=\Delta x$, $\Delta t=2\Delta x$ and $\Delta t=4\Delta x$. We compute the $L_2$ and $L_{\infty}$ errors with each neural network solver at final time $T=0.1$. Errors are listed in Table \ref{test4-table5}. The three network solvers give similar errors. The accuracy of CANN method seems to be independent of time step size $\Delta t$.

\begin{table}[!htb]
\centering
\begin{tabular}[c]{l l l l}
\hline
$\Delta t$ & $L_2$ & $L_{\infty}$ \\
\hline
4$\Delta x$ & $2.1981e^{-3}$ & $6.8272e^{-3}$\\
2$\Delta x$ & $2.4969e^{-3}$ & $7.2399e^{-3}$ \\
$\Delta x$ & $2.4736e^{-3}$ &$7.2650e^{-3}$ \\
\hline
\end{tabular}
\caption{Errors of neural network methods for Heat equation (Case II), $\Delta  x=1/160$ fixed, different $\Delta  t$}\label{test4-table5}
\end{table}

\vspace{.05in}
\noindent
{\bf Case III:}

Motivated by infinite speed of propagation, we wonder if the error and accuracy may be improved with increased stencil width of the network input vector, the $p$ and $q$ values of (\ref{nn:input-vector-stencil}). Three spatial mesh sizes of $\Delta x=1/40, 1/80, 1/160$ with $\Delta t=\Delta x$ are studied. Different to case I, we gradually increase the input vector stencil width with refined mesh. We have $p=q=2$ for $\Delta x=1/40$, $p=q=4$ for $\Delta x=1/80$ and $p=q=8$ for $\Delta x=1/160$. The three networks use same spatial cells to evolve the solution average. Errors are listed in Table \ref{test4-table6}. There is no sign of improvement with wider stencil included.
\begin{table}[!htb]
\centering
\begin{tabular}[c]{l l l l}
\hline
$\Delta x$ & $L_2$ & $L_{\infty}$ \\
\hline
1/40 & $7.1179e^{-3}$ & $1.8046e^{-2}$\\
1/80 & $6.3502e^{-3}$ & $1.8319e^{-2}$ \\
1/160 & $6.4138e^{-3}$ &$2.1510e^{-2}$ \\
\hline
\end{tabular}
\caption{Errors of neural network methods for Heat equation (Case III), $\Delta  t=\Delta x$, wider stencil on refined mesh}\label{test4-table6}
\end{table}
%%%%%%%%%%%%%%%%%%%%%%%%%%%%%%%%%%%%%%%%%%

\begin{example}\label{ex10} {\bf \emph {Linear convection diffusion equation}}
\end{example}
In this subsection, we consider linear convection-diffusion equation
\begin{equation*}
    u_t=u_{xx}+u_x, ~~~x\in D,~~t\geq 0,
\end{equation*}
with initial $u(x,0)=\sin(x)$. Spatial domain is $D=[0,2\pi]$. We further check whether the first order of accuracy can be obtained. We choose same network input vector of (\ref{eq:heat-input-vector}) as for Heat equation. Periodic boundary conditions are considered. Exact solution is available with  $u(x,t)=e^{-t}\sin(x+t)$. The picked network structure involves 1 hidden layer and 15 neurons. Number of iterations is taken as $K=5\times 10^6$. We adapt final time $T=\pi/4$ to compute all errors and orders. 

For the order of convergence test, four spatial mesh sizes of $\Delta x= \frac{\pi}{40}, \frac{\pi}{80}, \frac{\pi}{160}, \frac{\pi}{320}$ are studied. For each $\Delta x$, we choose $\Delta t=\Delta x$ correspondingly. Again, the four neural network solvers are similar to each other except the mesh size. In Table \ref{test10-table1}, we list the $L_2$ and $L_{\infty}$ errors. Roughly first order of convergence is observed.

\begin{table}[!htb]
\centering
\begin{tabular}[c]{l l l l l l}
\hline
$\Delta x$ & $L_2$ & order & $L_{\infty}$ & order\\
\hline
$\pi/40$  & $5.3013e^{-3}$ & & $3.2397e^{-3}$&\\
$\pi/80$ & $1.3801e^{-3}$ &1.94 & $7.5132e^{-4}$ &2.11\\
$\pi/160$ & $5.0091e^{-4}$ &1.46 & $2.6084e^{-4}$ &1.52\\
$\pi/320$ & $2.6771e^{-4}$ &0.91 & $1.5045e^{-4}$ &0.89\\
\hline
\end{tabular}
\caption{Errors and orders of neural network methods for linear convection diffusion equation, $\Delta t=\Delta x$.}\label{test10-table1}
\end{table}

The second test is similar to the Case II of Example \ref{ex4}. We fix the spatial size $\Delta x=\pi/160$ and vary the time step size from $\Delta t=\Delta x$, $\Delta t=2\Delta x$ to $\Delta t=4\Delta x$. The computed $L_2$ and $L_{\infty}$ errors of the three network solvers are listed in Table \ref{test10-table2}. The errors seem to be related to time step size $\Delta t$, even the relationship is not clear. Again, the cell-average based neural network allows large time step. Notice the choice of $\Delta t=4\Delta x$ involves a time step size roughly $200$ times bigger than the regular CFL restriction of $\Delta t\approx \Delta x^2$.

\begin{table}[!htb]
\centering
\begin{tabular}[c]{l l l l}
\hline
$\Delta t$ & $L_2$ & $L_{\infty}$ \\
\hline
4$\Delta x$ & $1.0249e^{-3}$ & $5.8832e^{-4}$\\
2$\Delta x$ & $8.2225e^{-4}$ & $4.2640e^{-4}$ \\
$\Delta x$ & $5.0091e^{-4}$ &$2.6084e^{-4}$ \\
$\Delta x/2$ & $3.7288e^{-4}$ &$2.0151e^{-4}$ \\
\hline
\end{tabular}
\caption{Errors of neural network methods for linear convection diffusion equation, $\Delta  x=\pi/160$ fixed, different $\Delta  t$}\label{test10-table2}
\end{table}

% \begin{figure}[!htb]\label{test-4}
% \begin{minipage}[c]{0.50\linewidth}
% \centering
% \includegraphics[width=1\linewidth]{fig/convection-diffusion.eps}
% \caption{Trajectory approximations for different long time run}\label{test4-figure}
% \end{minipage}
% \begin{minipage}[c]{0.40\linewidth}
% \centering
% \begin{tabular}[c]{l l l l l}
% \hline
%   & $t=0.133$ & $t=0.3$ & $t=1$ \\
% \hline
% $L_2$  & $2.9700e^{-6}$ & $5.0816e^{-6}$ & $1.2069e^{-5}$\\
% \hline
% $L_{\infty}$ & $7.3276e^{-5}$ & $1.1315e^{-4}$ & $2.1859e^{-4}$\\
% \hline
% \end{tabular}
% \caption{$L_2$ error and $L_{\infty}$ error for different long time simulations}\label{test4-table}
% \end{minipage}
% \end{figure}
%%%%%%%%%%%%%%%%%%%%%%%%%%%%%%%%%%%
\subsection{Nonlinear hyperbolic equation}\label{S3-3}
In this section, we investigate the effectiveness of neural network method solving nonlinear conservation laws. We consider the inviscid Burgers’ equation 
\begin{equation}
    u_t+\left(\frac{u^2}{2}\right)_x=0, ~~(x, t)\in D\times R^{+}.\label{ex-inviscid-burgers}
\end{equation}
Four benchmark problems are studied to illustrate the computational challenges of nonlinear hyperbolic PDEs. We first consider the smooth initial $\sin(x)$ wave which develops into a shock after finite time evolution. Then we study two Riemann problems with two piece wise constants initials that will develop into a shock and a rarefaction wave. Last example involves three piece wise constants initial that will develop into the interaction between rarefaction wave and shock. In all four examples, Dirichlet boundary conditions are applied.

For nonlinear problems, multiple time levels of data pairs are necessary and applied training the neural network. For all four examples, time step size $\Delta t=0.1$ is taken. We roughly have $\Delta t\approx 2\Delta x$. Solution cell average values up to $t=2.0$, as listed below
\begin{equation}\label{data-pair-burgers}
S =\left\{\left(\bar{u}_{j}^n,\bar{u}_{j}^{n+1}\right),\,\,j=1, \cdots, J\right\}_{n=0}^{m=19},
\end{equation}
are included in the training data set. So we have twenty time levels of solution average pairs $\left(\bar{u}_{j}^n,\bar{u}_{j}^{n+1}\right)$ used in training the network. Spatial domain is either $D=[0, 2\pi]$ or $D=[-1, 5]$. Mesh size of $\Delta x=\frac{2\pi}{100}$ for example \ref{ex5} and $\Delta x=\frac{6}{100}$ for other three examples are applied. 

\begin{figure}[htbp]
\centering
\subfigure[Example \ref{ex5}: sine initial]{
\includegraphics[width=0.35\linewidth]{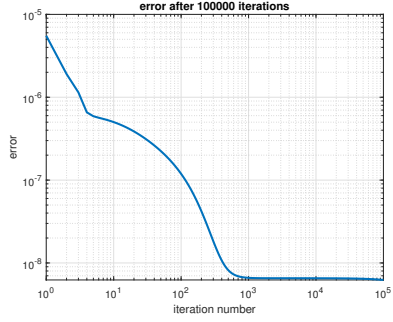}
%\caption{}
}\quad\quad
\subfigure[Example \ref{ex6}: shock wave]{
\includegraphics[width=0.35\linewidth]{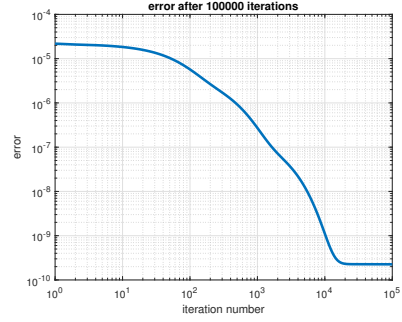}
%\caption{}
}
\quad
\subfigure[Example \ref{ex7}: rarefaction wave]{
\includegraphics[width=0.35\linewidth]{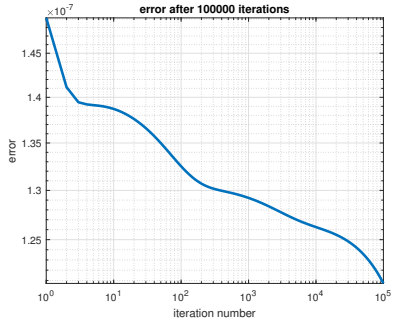}
%\caption{}
}\quad\quad
\subfigure[Example \ref{ex8}: rarefaction and shock interaction]{
\includegraphics[width=0.35\linewidth]{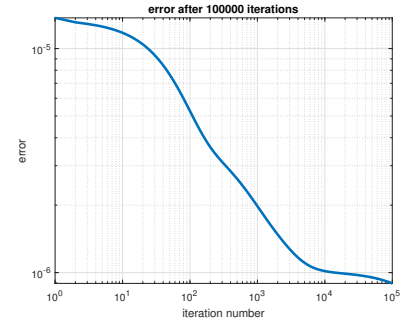}
%\caption{}
}
\caption{Invisid Burgers' equation neural network squared $L_2$ training errors (\ref{nn:L2-error-square}) at the last time level
} 
\label{fig:inviscid-Burgers-error}
\end{figure}

%%%%%%%%%%%%%%%%%%%%%%%%%%%%%%%%%%%%%%%%%%%%%%%%%%%%%%%%%%%%%%%%%%%%%%%%
\begin{example}\label{ex5} {\bf \emph {Smooth initial sine wave}}
\end{example}
We first consider the inviscid Burgers equation (\ref{ex-inviscid-burgers}) associated with smooth initial value of  
\begin{equation*}
    u(x, 0)=\sin(x),
\end{equation*}
over domain $D=[0, 2\pi]$. Zero boundary conditions $u(0, t)=u(2\pi, t)=0$ and its zero extension to out of domain ghost cells are applied. With $\Delta t=0.1$ and $\Delta x=2\pi/100$ and characteristic speed less than one, the network input vector is taken as 
$$\vv_{j}^{in} = \Big[\uAve^n_{j-3}, \uAve^n_{j-2}, \uAve^n_{j-1}, \uAve^n_j, \uAve^n_{j+1}, \uAve^n_{j+2},\uAve^n_{j+3}\Big]^T.$$

\begin{figure}[!htb]
\centering
\includegraphics[width=0.3\linewidth]{fig/inviscid_Burgers_sin/u_0.eps}\quad\quad
\includegraphics[width=0.3\linewidth]{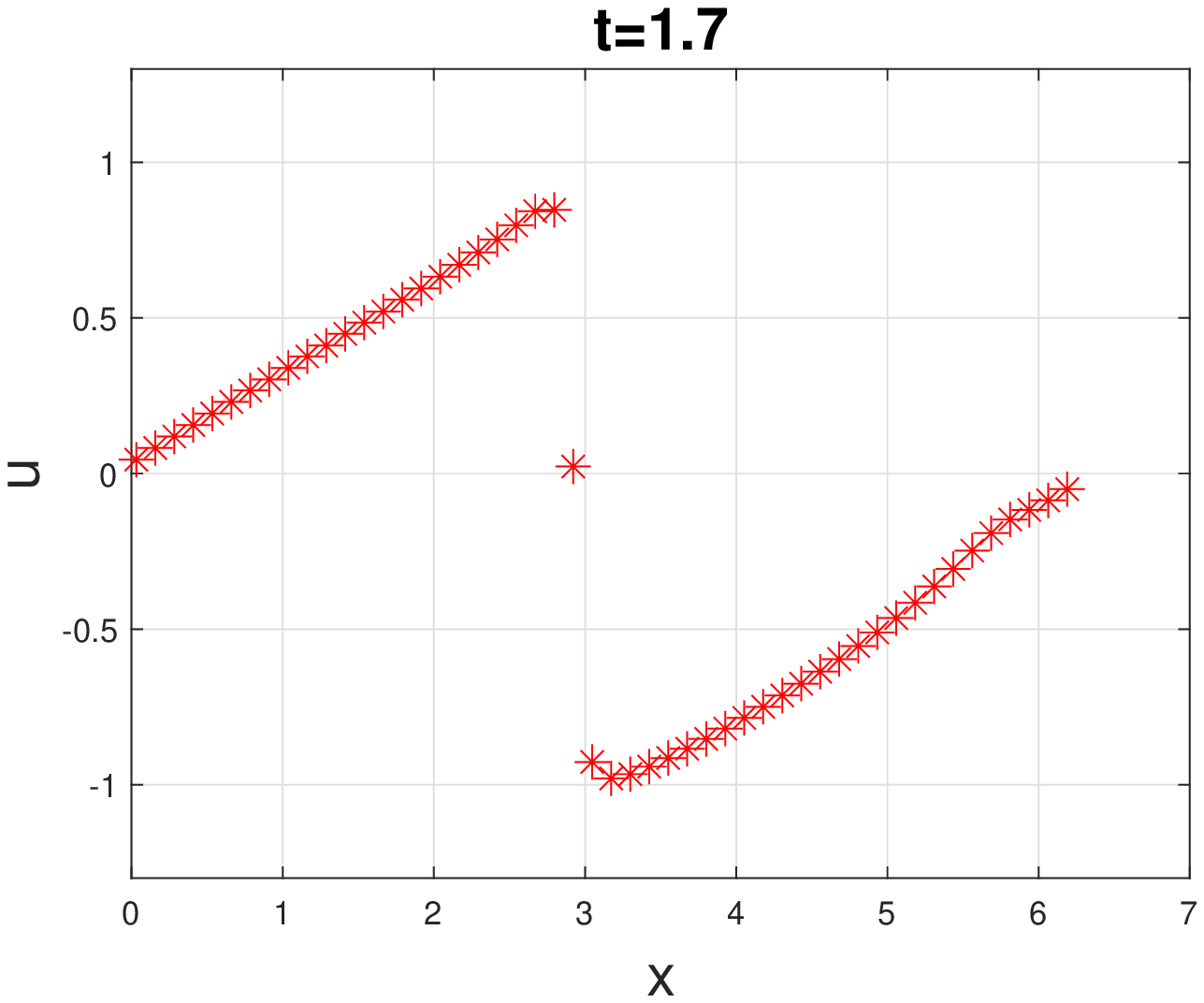}\\
\includegraphics[width=0.3\linewidth]{fig/inviscid_Burgers_sin/u26.eps}\quad\quad
\includegraphics[width=0.3\linewidth]{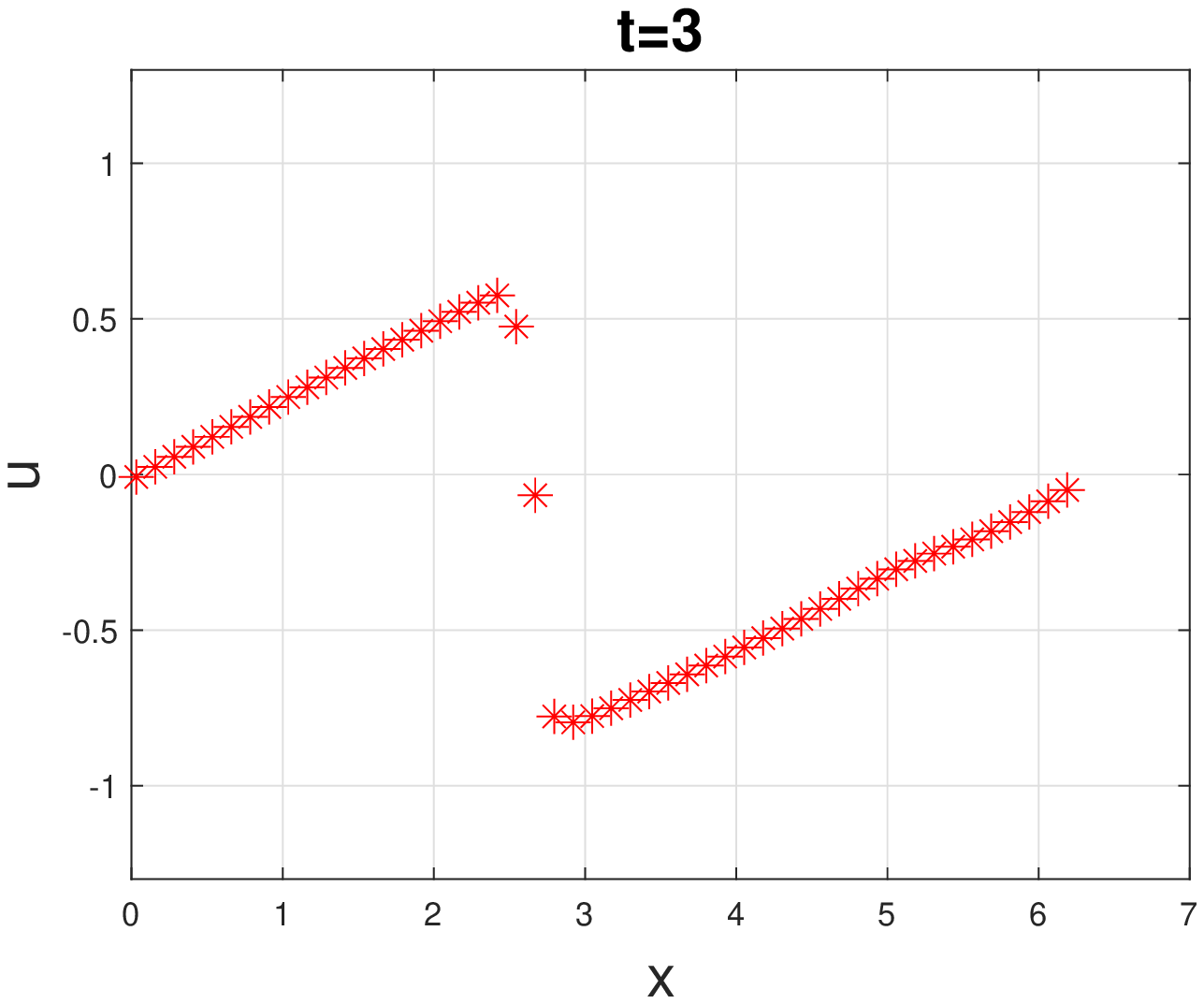}
\caption{neural network simulation for smooth $\sin(x)$ wave initial (Example \ref{ex5})}
\label{fig:inviscid-Burgers-sin}
\end{figure}

This stencil includes the characteristic and the domain of dependence. Training data pairs $\left(\bar{u}_{j}^n,\bar{u}_{j}^{n+1}\right)$ are obtained from highly accurate discontinuous Galerkin method. The network itself contains 2 hidden layers with 8 neurons per layer. The network training is conducted for up to $K=10^5$ iterations. We also output the squared $L_2$ training error of (\ref{nn:L2-error-square}) corresponding to iteration at the very last time level $t=2$, see Fig \ref{fig:inviscid-Burgers-error}.

Well trained network is applied to simulate solution evolution to $T=3$. Before and after shock developed are illustrated in Fig \ref{fig:inviscid-Burgers-sin}. The cell-average based neural network results are comparable to those obtained with classical numerical methods. Shock evolution is sharply captured and no oscillation is generated with the neural network method.

%%%%%%%%%%%%%%%%%%%%%%%%%%%%%%%%%%%%%%%%%%%%%%%%%%%%%%%%%%%%%%%%%%%

\begin{example}\label{ex6} {\bf \emph {Single shock propagation}}
\end{example}
In this example, we consider a Riemann problem of (\ref{ex-inviscid-burgers}) with piece wise constant initial 
\begin{equation*}
    u(x, 0)=\left\{
             \begin{array}{lr}
             1,  ~~~~x<0,\\
             0,  ~~~ otherwise. 
             \end{array}
\right.
\end{equation*}
Computational domain is $D=[-1, 5]$. Dirichlet boundary conditions of $u(-1, t)=1$, $u(5, t)=0$ and its out of domain extensions are applied. Cell size $\Delta x=0.06$ and time step size $\Delta t=0.1$ are taken. The network input vector of (\ref{nn:input-vector-stencil}) is chosen as 
$$\vv_{j}^{in} = \Big[\uAve^n_{j-4}, \uAve^n_{j-3}, \uAve^n_{j-2}, \uAve^n_{j-1}, \uAve^n_j, \uAve^n_{j+1}, \uAve^n_{j+2}\Big]^T.$$
This is a biased choice which includes more points to the left. Again the choice covers the domain of dependence. Training data are obtained from the exact solution 
\begin{equation}\label{exact_solution_shock}
        u(x, t)=\left\{
             \begin{array}{lr}
             1,  ~~~x\leq \frac{1}{2}t,\\
             0,  ~~~ otherwise. 
             \end{array}
\right.
\end{equation}

\begin{figure}[!htb]
\centering
\includegraphics[width=0.3\linewidth]{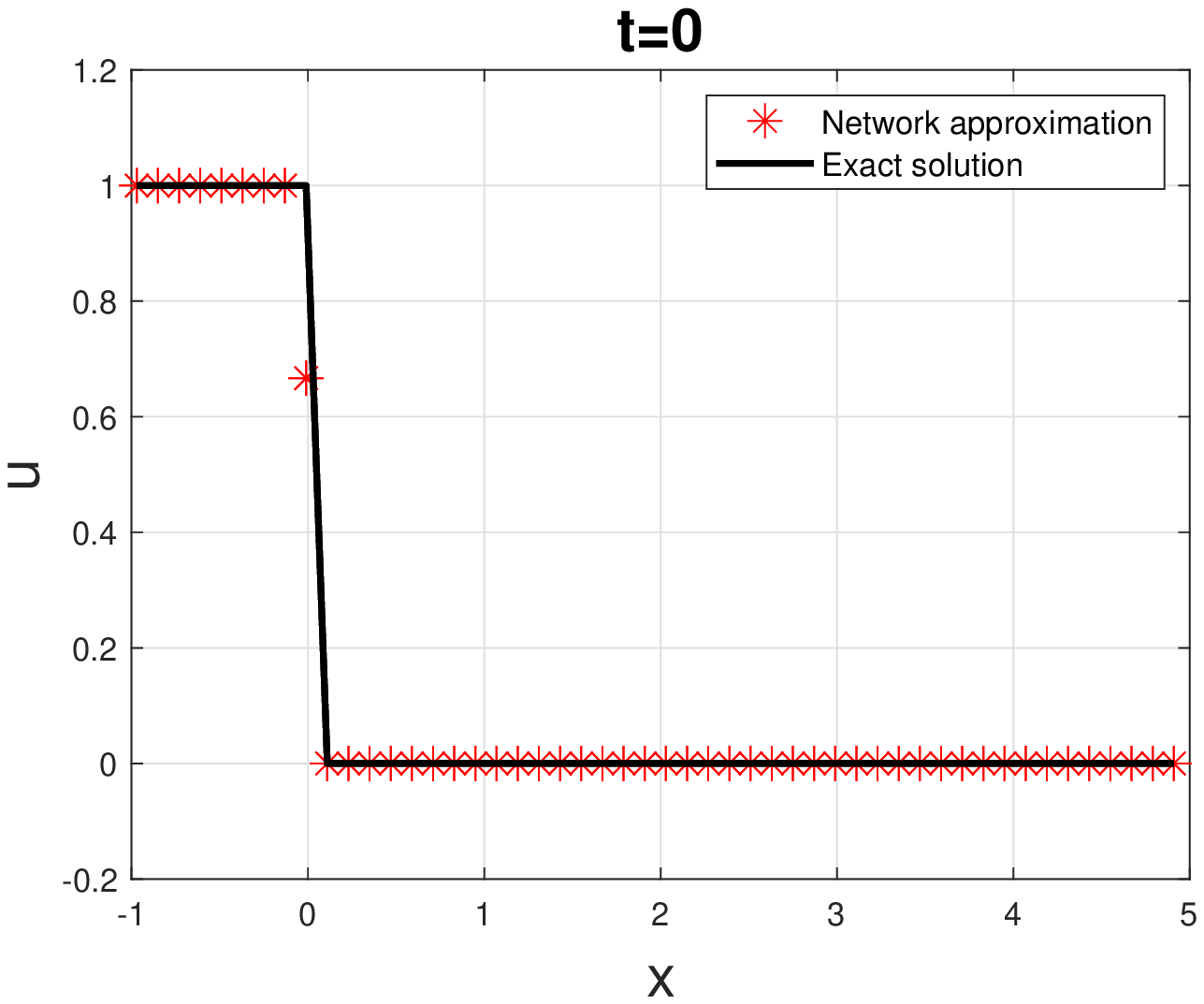}\quad\quad
\includegraphics[width=0.3\linewidth]{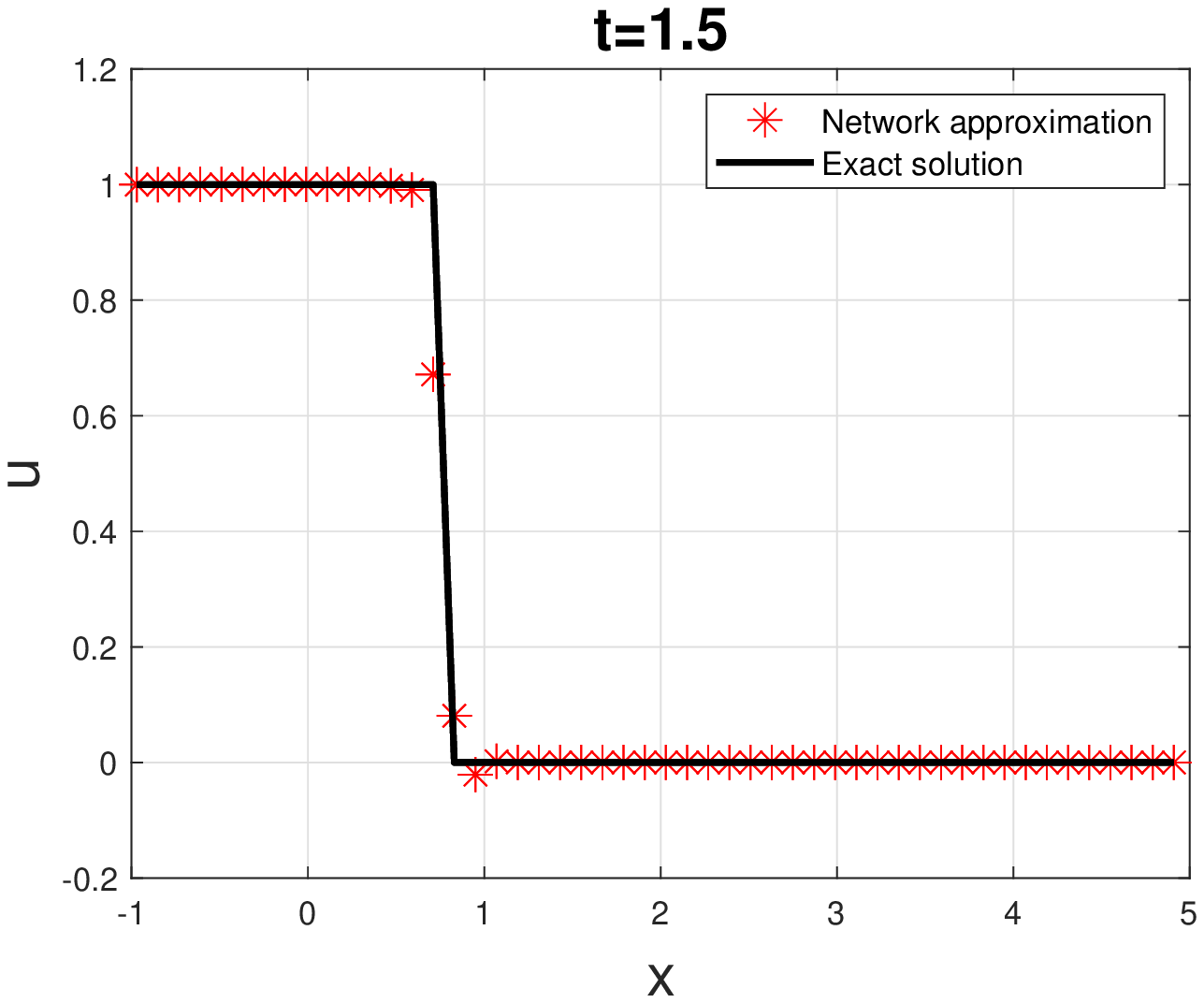}\\
\includegraphics[width=0.3\linewidth]{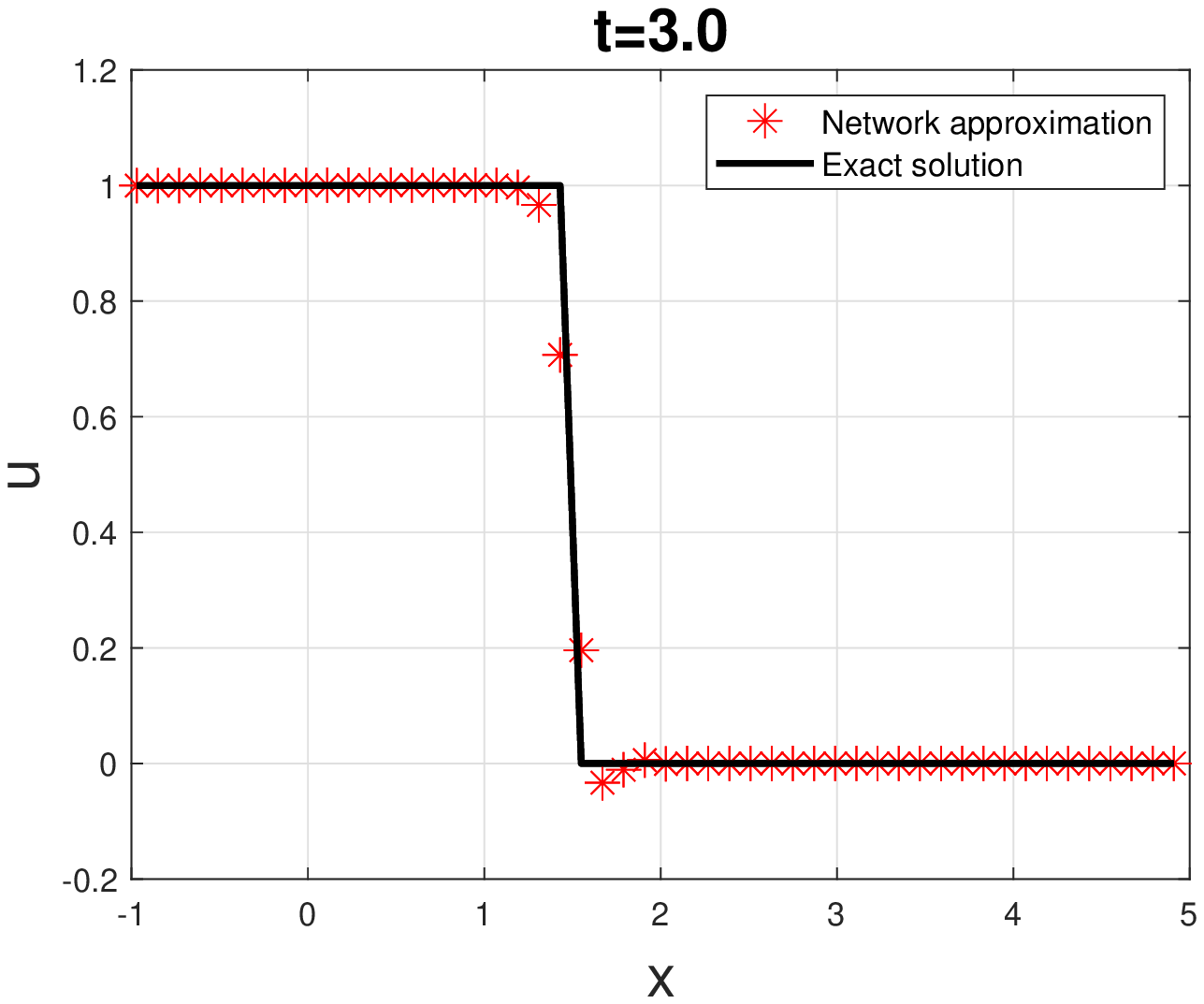}\quad\quad
\includegraphics[width=0.3\linewidth]{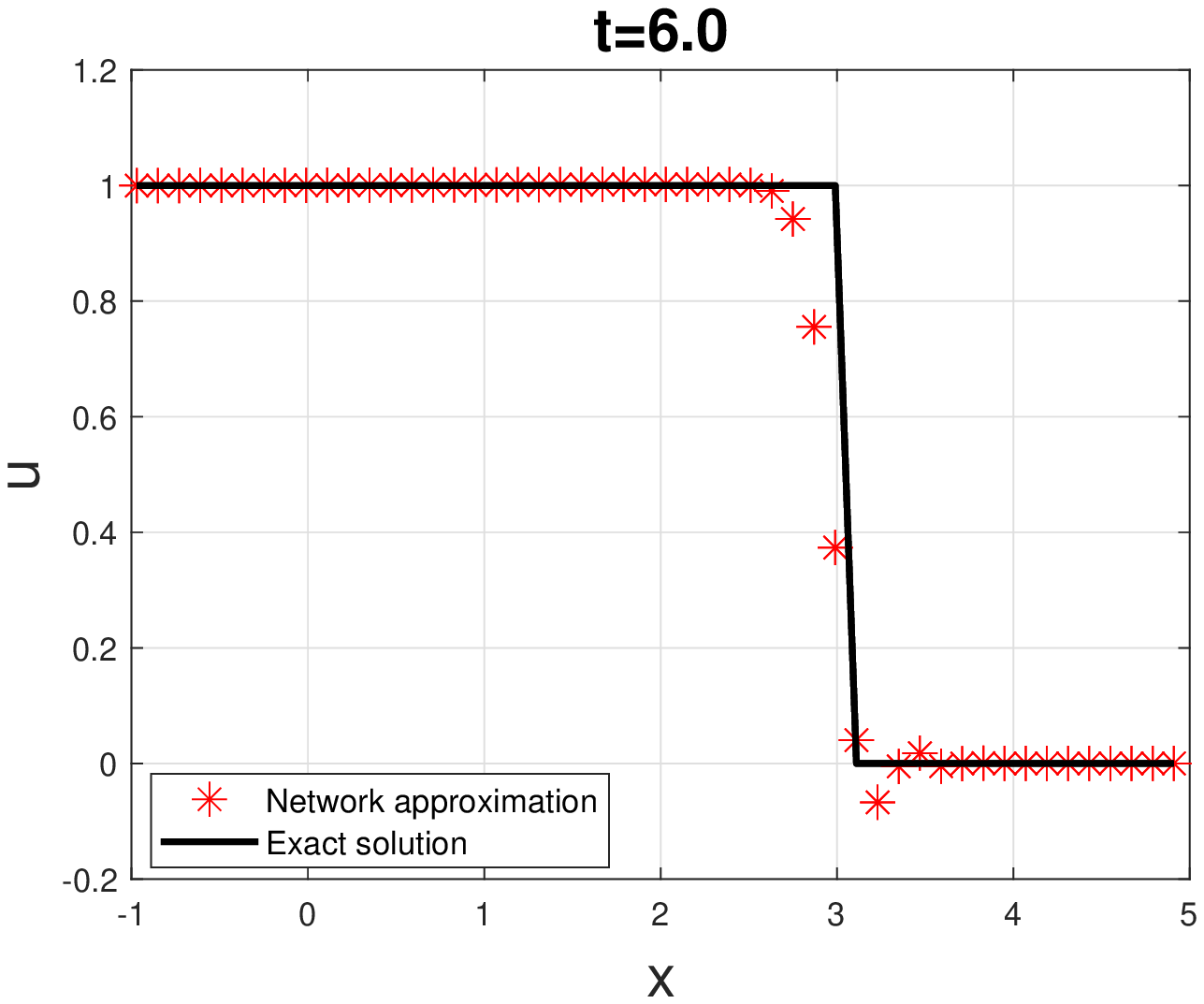}
\caption{shock propagation (Example \ref{ex6}) with neural network method}
\label{fig:inviscid-Burgers-shock}
\end{figure}
The neural network picked contains 1 hidden layer of 8 neurons. Network training is conducted for up to $K=10^5$ iterations. The squared $L_2$ training error of (\ref{nn:L2-error-square}) to iteration at $t=2$ are output in Fig \ref{fig:inviscid-Burgers-error}. After well trained, the neural network is applied solving the Riemann problem to $T=8.0$. Screen shots of $t=0$, $t=1.5$, $t=3.0$ and $t=6.0$ are shown in Figure \ref{fig:inviscid-Burgers-shock}. The neural network method can accurately and sharply capture the shock evolution. Notice the neural network simulation at $t=8$ is way later than the training data set of $t=2$.
%%%%%%%%%%%%%%%%%%%%%%%%%%%%%%%%%%%%%%%%%%%%%%%%%%%%

\begin{example}\label{ex7} {\bf \emph {Rarefaction wave}}
\end{example}
In this example, we consider a Riemann problem of (\ref{ex-inviscid-burgers}) with piece wise constant initial 
\begin{equation*}
    u(x, 0)=\left\{
             \begin{array}{lr}
             0,  ~~~x<0,\\
             1,  ~~~ otherwise. 
             \end{array}
\right.
\end{equation*}
Domain is $D=[-1, 5]$. Dirichlet boundary condition $u(-1, t)=0$, $u(5, t)=1$ and its out of domain extension are applied. Cell size $\Delta x=0.06$ and time step size $\Delta t=0.1$ are chosen. Network input vector is taken as
$$\vv_{j}^{in} = \Big[ \uAve^n_{j-2}, \uAve^n_{j-1}, \uAve^n_j, \uAve^n_{j+1}\Big]^T.$$ 
The choice still includes the domain of dependence. Training data are generated from exact solution
\begin{equation*}
%\label{exact_solution_rarefaction}
        u(x, t)=\left\{
             \begin{array}{lr}
             0,  ~~x<0,\\
             \frac{x}{t},  ~~0\le x\le t,\\
             1,  ~~~ otherwise. 
             \end{array}
\right.
\end{equation*}
\begin{figure}[!htb]
\centering
\includegraphics[width=0.3\linewidth]{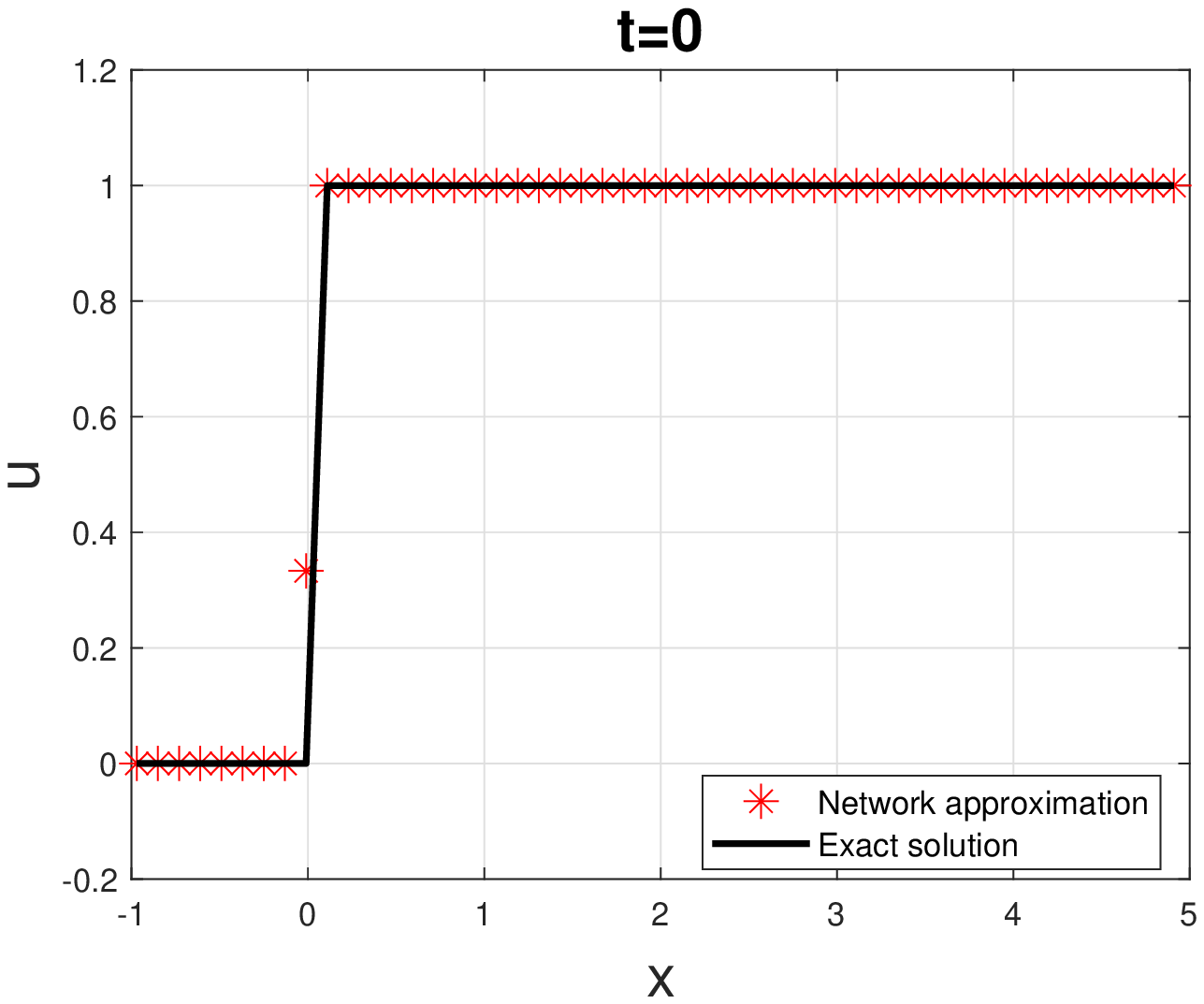}\quad\quad
\includegraphics[width=0.3\linewidth]{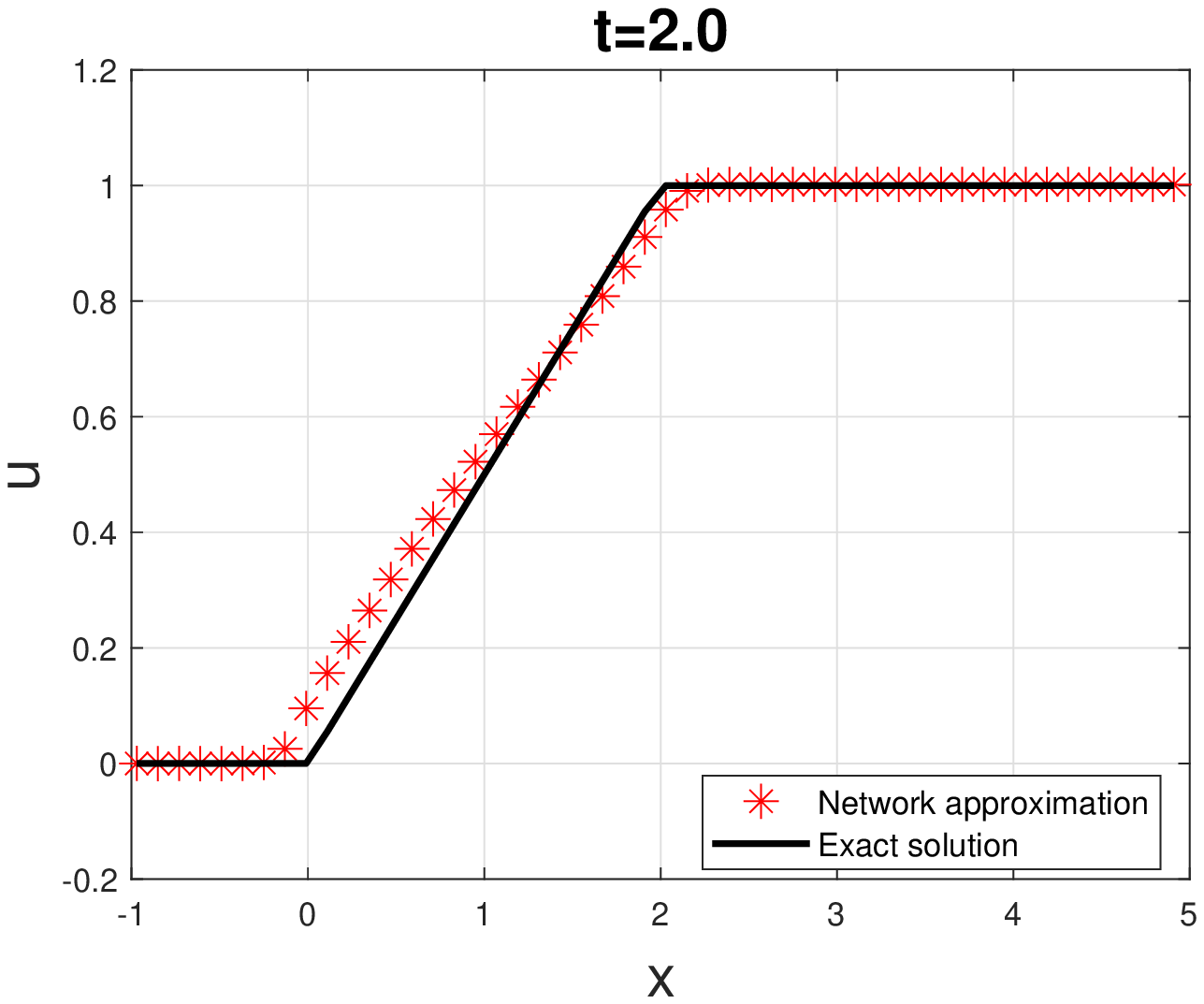}\\
\includegraphics[width=0.3\linewidth]{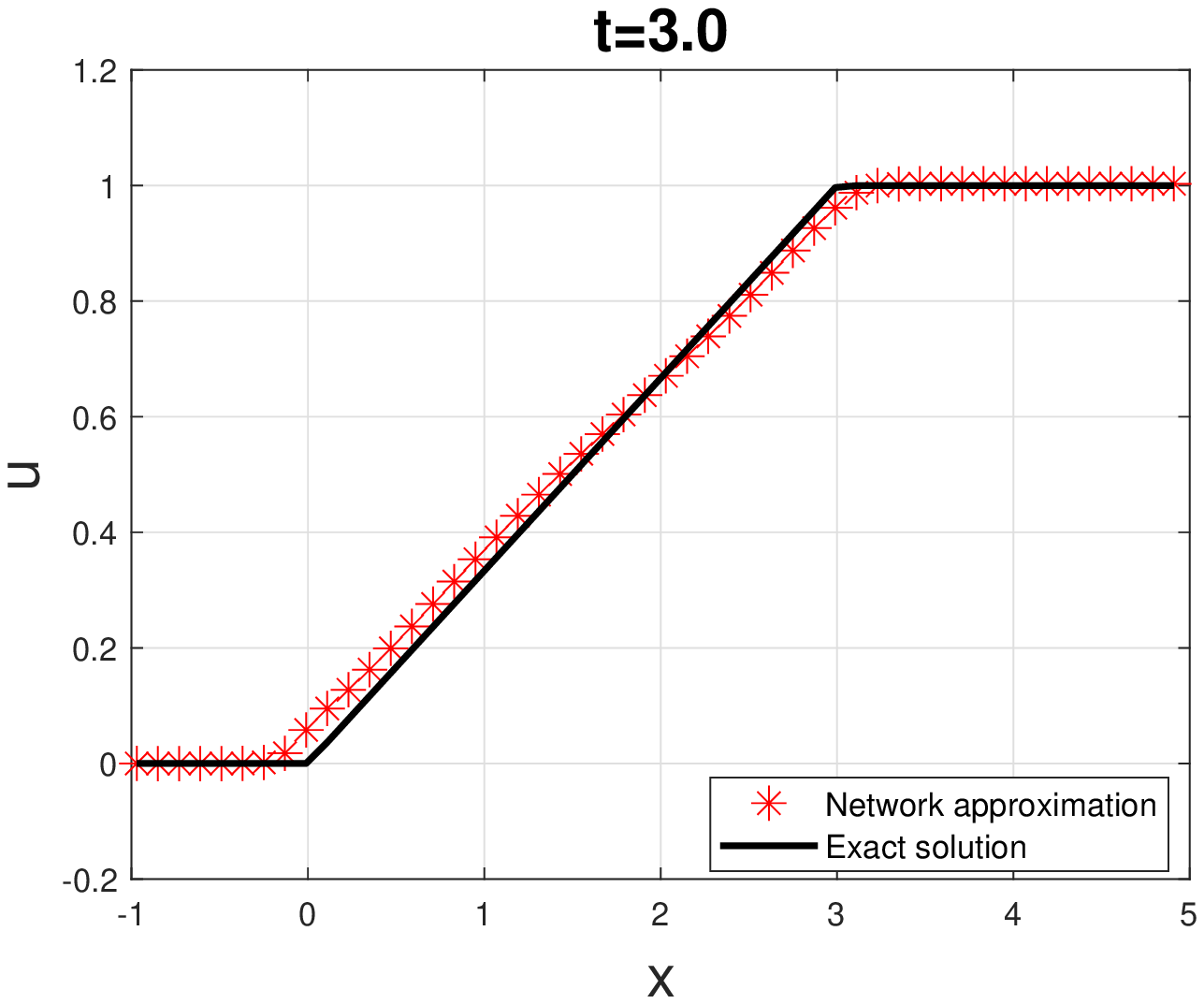}\quad\quad
\includegraphics[width=0.3\linewidth]{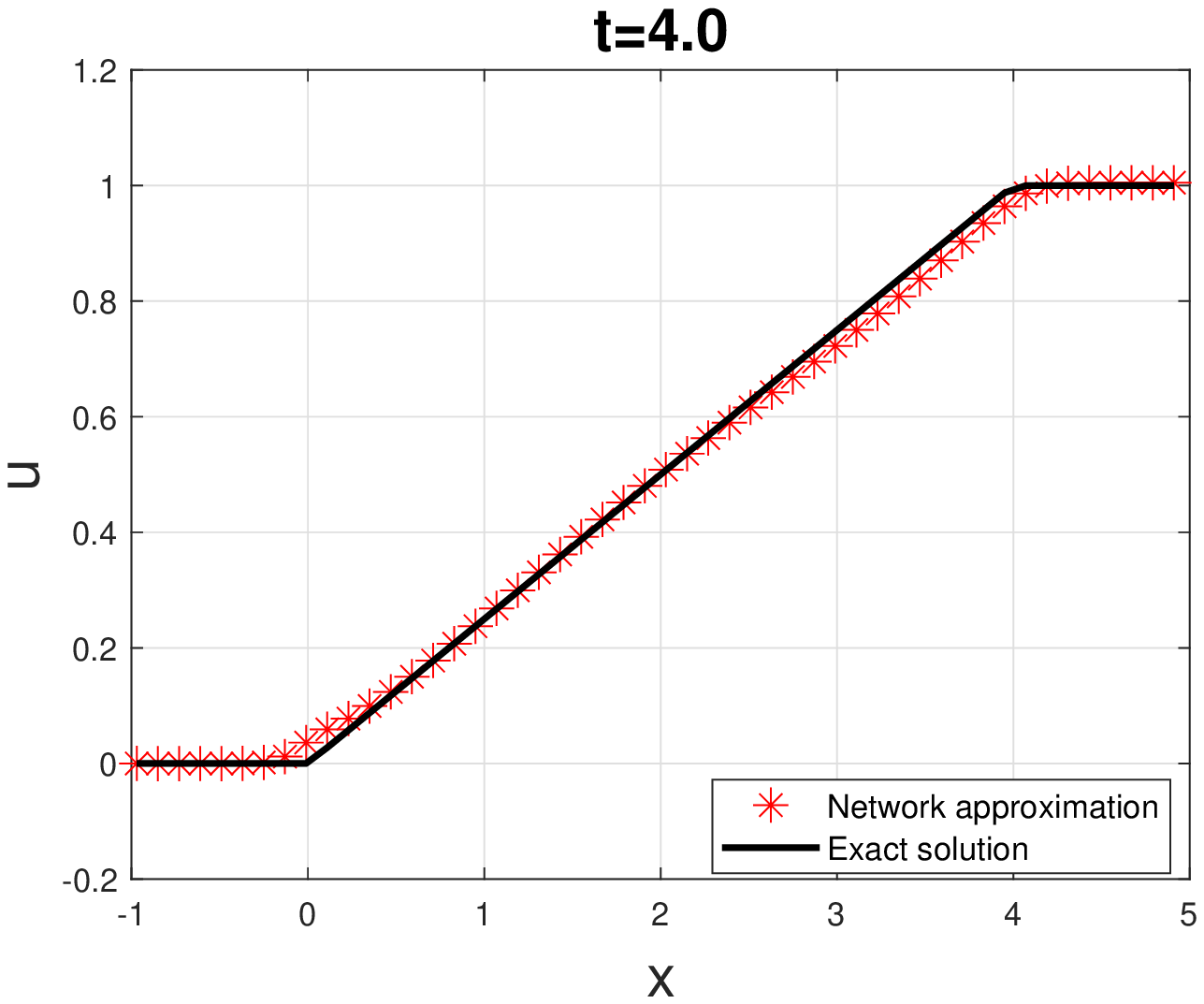}
\caption{rarefaction wave propagation (Example \ref{ex7}) with neural network method}
\label{fig:inviscid-Burgers-rarefaction}
\end{figure}

The network contains 2 hidden layers with 8 neurons per layer. Iteration steps up to $K=10^5$ are applied. The squared $L_2$ training error of (\ref{nn:L2-error-square}) corresponding to iteration at $t=2$ are also included in Figure \ref{fig:inviscid-Burgers-error}. After well trained, the neural network solver is applied solving this Riemann problem up to $T=4.0$. Screen shots are shown in Figure \ref{fig:inviscid-Burgers-rarefaction}. 
The neural network method can open up and well resolves the rarefaction wave evolution.

%%%%%%%%%%%%%%%%%%%%%%%%%%%%%%%%%%%%%%%%%%%%%%%%%%%%%%%%%%%%%%%%%%%%%%
\begin{example}\label{ex8} {\bf \emph {Interaction of rarefaction and shock waves}}
\end{example}
In this example, we consider solving (\ref{ex-inviscid-burgers}) with three piece wise constants initial 
\begin{equation*}
    u(x, 0)=\left\{
             \begin{array}{lr}
             0,  ~~~x<0,\\
             1,  ~~~0\leq x\leq 1,\\
             0,  ~~~ otherwise. 
             \end{array}
\right.
\end{equation*}
Spatial domain is set as $D=[-1, 5]$. Dirichlet boundary conditions $u(-1, t)=u(5, t)=0$ and its out of domain extension are applied. Network input vector is taken as 
$$\vv_{j}^{in} = \Big[\uAve^n_{j-4}, \uAve^n_{j-3}, \uAve^n_{j-2}, \uAve^n_{j-1}, \uAve^n_j, \uAve^n_{j+1}, \uAve^n_{j+2}\Big]^T.$$
Training data pairs $\left(\uAve_j^n,\uAve_j^{n+1}\right)$ are generated from the exact solution. Here the rarefaction wave travels faster than the shock and it merges into the shock wave after a while. Before the meet at $t=2$, the exact solution is given as
\begin{equation*}\label{exact:rarefaction_shock_1}
        u(x, t)=\left\{
             \begin{array}{lr}
             0,  ~~~x<0,\\
             \frac{x}{t},  ~~~0\le x<t,\\
             1,  ~~~t\le x\le 1+\frac{t}{2},\\
             0,  ~~~ otherwise.
             \end{array}
\right.
\end{equation*}
After $t > 2$, the rarefaction wave runs into the shock wave. The solution is composed with zero state to the left, followed by a rarefaction wave and connected with a shock to the right with which the shock speed slows down as time evolves. The shock speed is determined by the following formula
\begin{equation*}
    \sigma^{\prime}(t)=\frac{\frac{1}{2}(\frac{\sigma (t)}{t})^2}{\frac{\sigma (t)}{t}}=\frac{\sigma (t)}{2t}.
\end{equation*}
With initial value $\sigma (2)=2$, we obtain $\sigma (t)=\sqrt{2t}$. 
For $t>2$, the exact solution is given by
\begin{equation*}\label{exact:rarefaction_shock_2}
        u(x, t)=\left\{
             \begin{array}{lr}
             0,  ~~~x<0,\\
             \frac{x}{t},  ~~~0\le x\le\sqrt{2t},\\
             0,  ~~~ \sqrt{2t}<x. 
             \end{array}
\right.
\end{equation*}
\begin{figure}[htbp]
\centering
\includegraphics[width=0.3\linewidth]{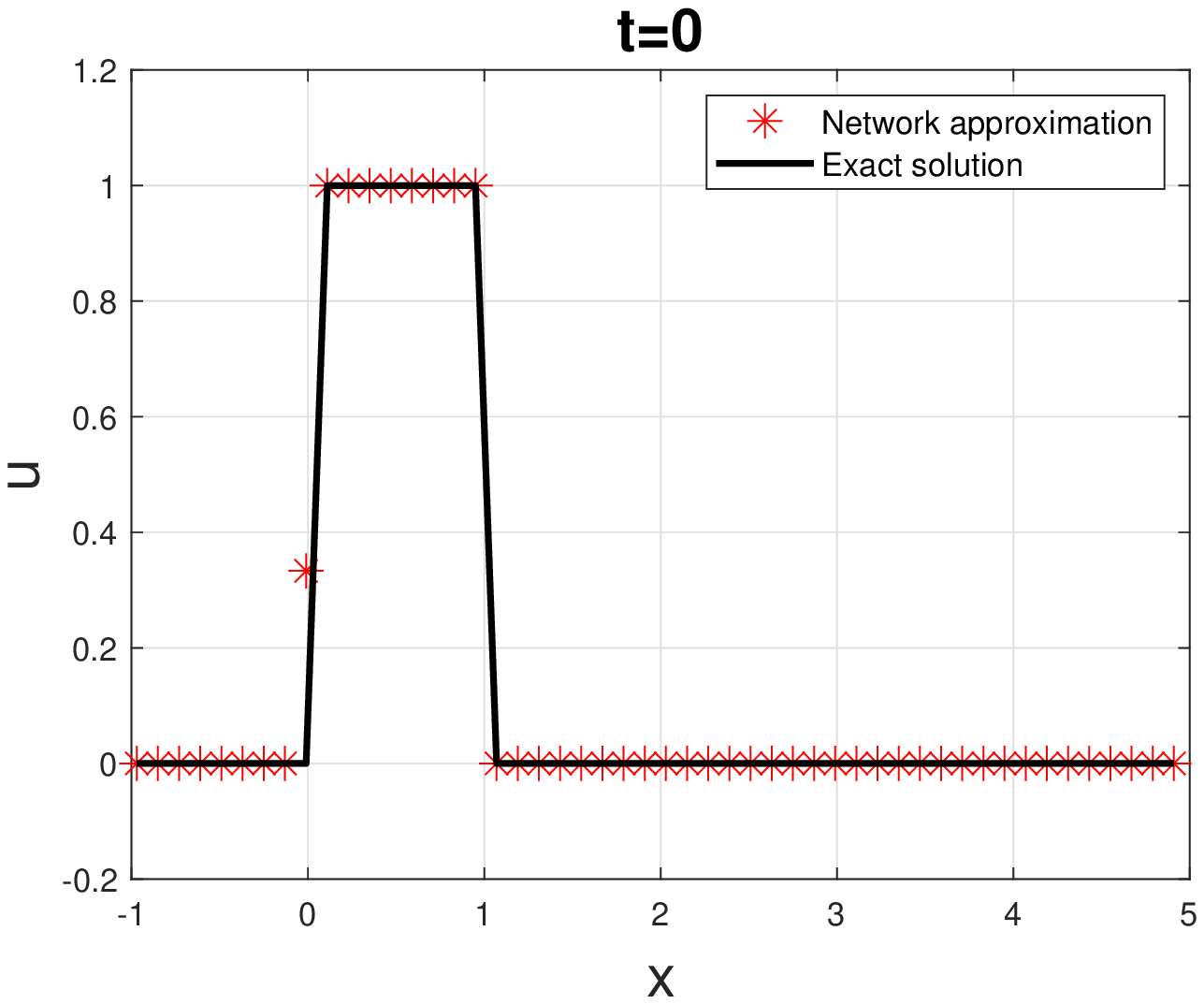}\quad\quad
\includegraphics[width=0.3\linewidth]{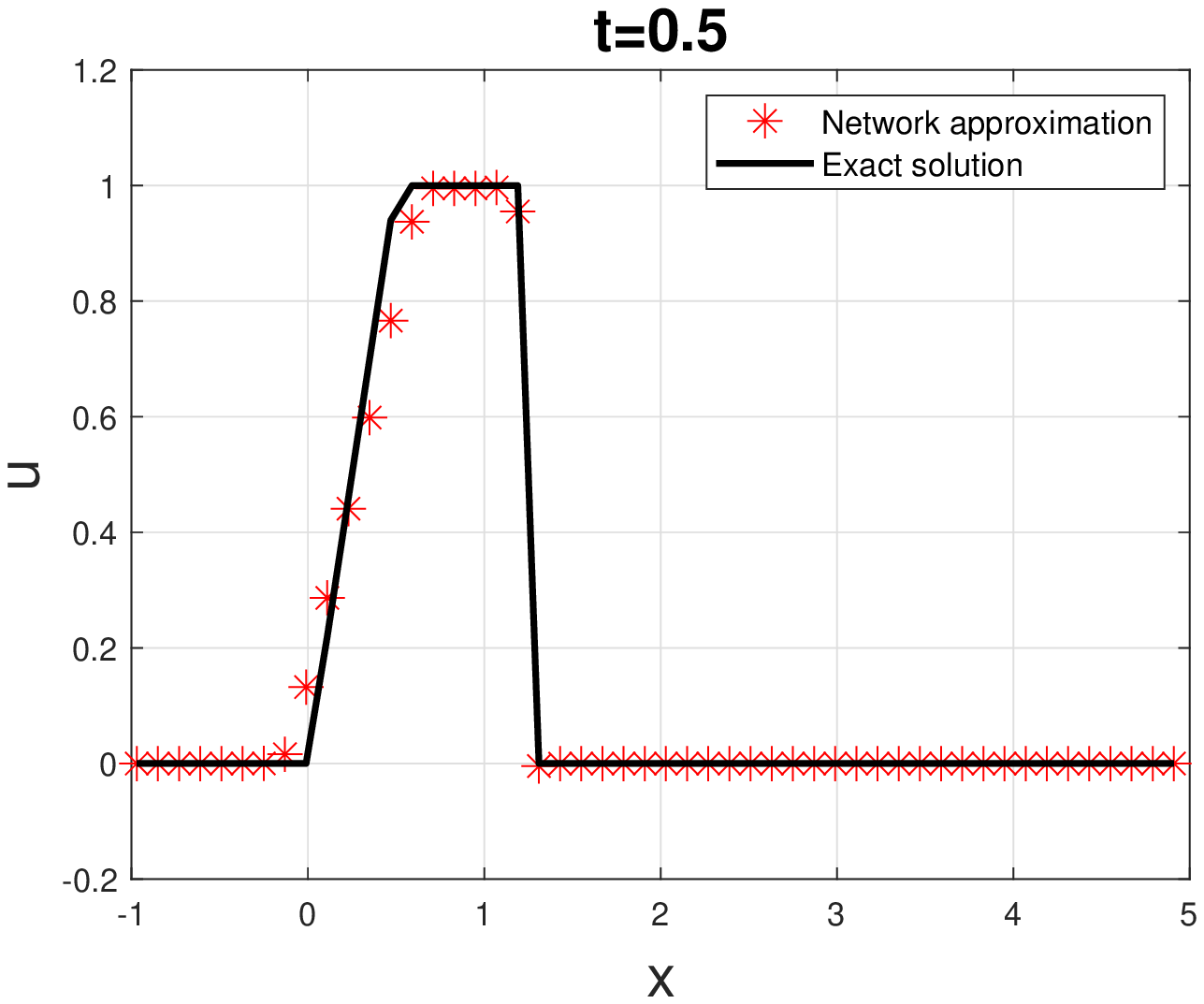}\\
\includegraphics[width=0.3\linewidth]{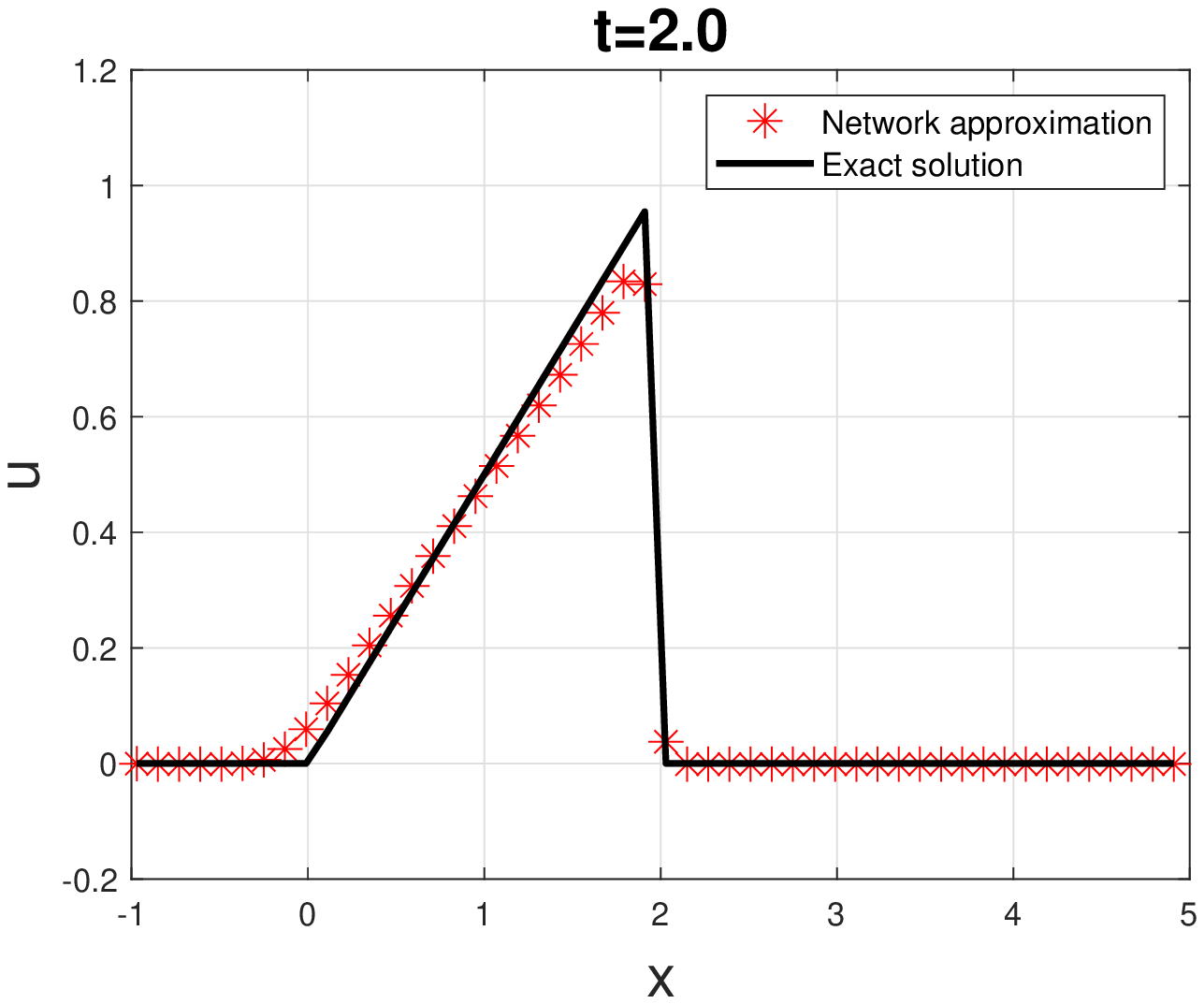}\quad\quad
\includegraphics[width=0.3\linewidth]{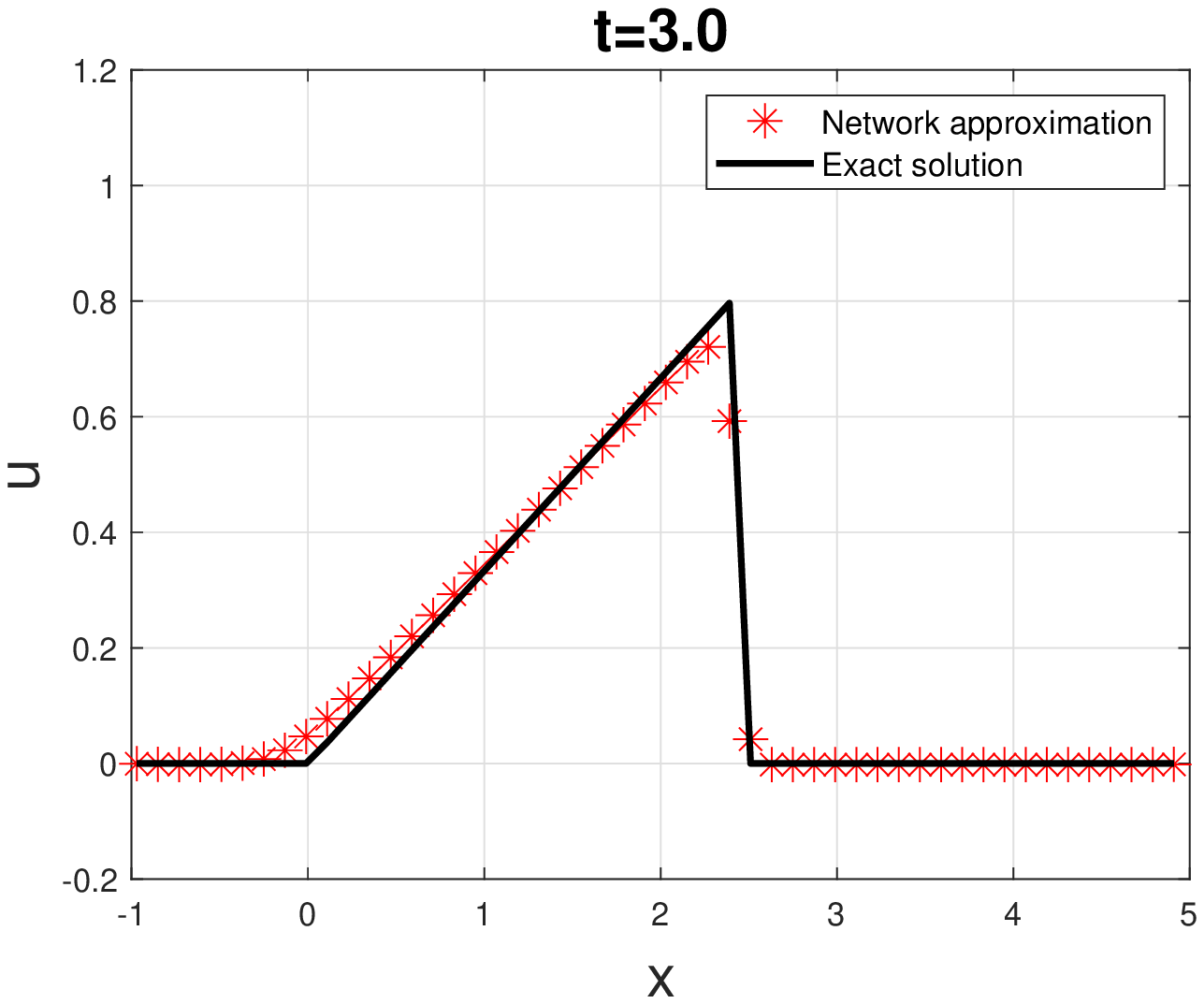}
\caption{rarefaction wave and shock wave interaction (Example \ref{ex8}) with neural network method}
\label{fig:inviscid-Burgers-rarefaction-shock}
\end{figure}

The picked network contains 2 hidden layers with 8 neurons per layer. Training is conducted up to $K=10^5$ iterations. Squared $L_2$ training error of (\ref{nn:L2-error-square}) corresponding to iteration at $t=2$ is listed in Figure \ref{fig:inviscid-Burgers-error}.
After well trained, neural network method is applied with its solution screen shots before and after waves interaction shown in Figure \ref{fig:inviscid-Burgers-rarefaction-shock}. Recall the training set only include solution pairs up to $t=2$ for which the rarefaction wave has not met the shock wave yet. The neural network method is capable of accurately capturing shock propagation after the interaction. 
We mention that we also check out the $L_2$ errors, which are around $O(10^{-3})$ for all four examples, at a time when solutions all develop singularities. Recall the cell size is around $\Delta x\approx 0.06$. 

\subsection{Nonlinear convection diffusion equation}\label{S3-5}

In this section, we investigate the viscous Burgers’ equation 
\begin{equation}
    u_t+\left(\frac{u^2}{2}\right)_x=\mu u_{xx}, ~~(x, t)\in D\times R^{+},\label{ex-viscous-burgers}
\end{equation}
with $\mu=0.1$. Zero boundary conditions of $u(0, t)=u(2\pi, t)=0$ are applied. Same time step size $\Delta t= 0.1$ and cell size $\Delta x=\frac{2\pi}{100}$ as in Example \ref{ex5} are taken. The cell-average based neural network method can adapt large time step size $\Delta t$. We use this setting to simply illustrate the effectiveness and efficiency of neural network method solving nonlinear convection diffusion equations. The network input vector is taken the same as that of Example \ref{ex5} 
$$\vv_{j}^{in} = \Big[\uAve^n_{j-3}, \uAve^n_{j-2}, \uAve^n_{j-1}, \uAve^n_j, \uAve^n_{j+1}, \uAve^n_{j+2},\uAve^n_{j+3}\Big]^T.$$ 
Training data $\left\{\uAve_j^{n},\uAve_j^{n+1}\right\}$ are generated from the refined mesh highly accurate discontinuous Galerkin method. With time step $\Delta t=0.1$ and up to $t=2.0$, twenty time levels of solution averages are applied in the training data set. The network structure contains 2 hidden layers of each with 8 neurons. The network training is conducted for up to $K=10^5$ iterations.

\begin{figure}[htbp]
\centering
\includegraphics[width=0.3\linewidth]{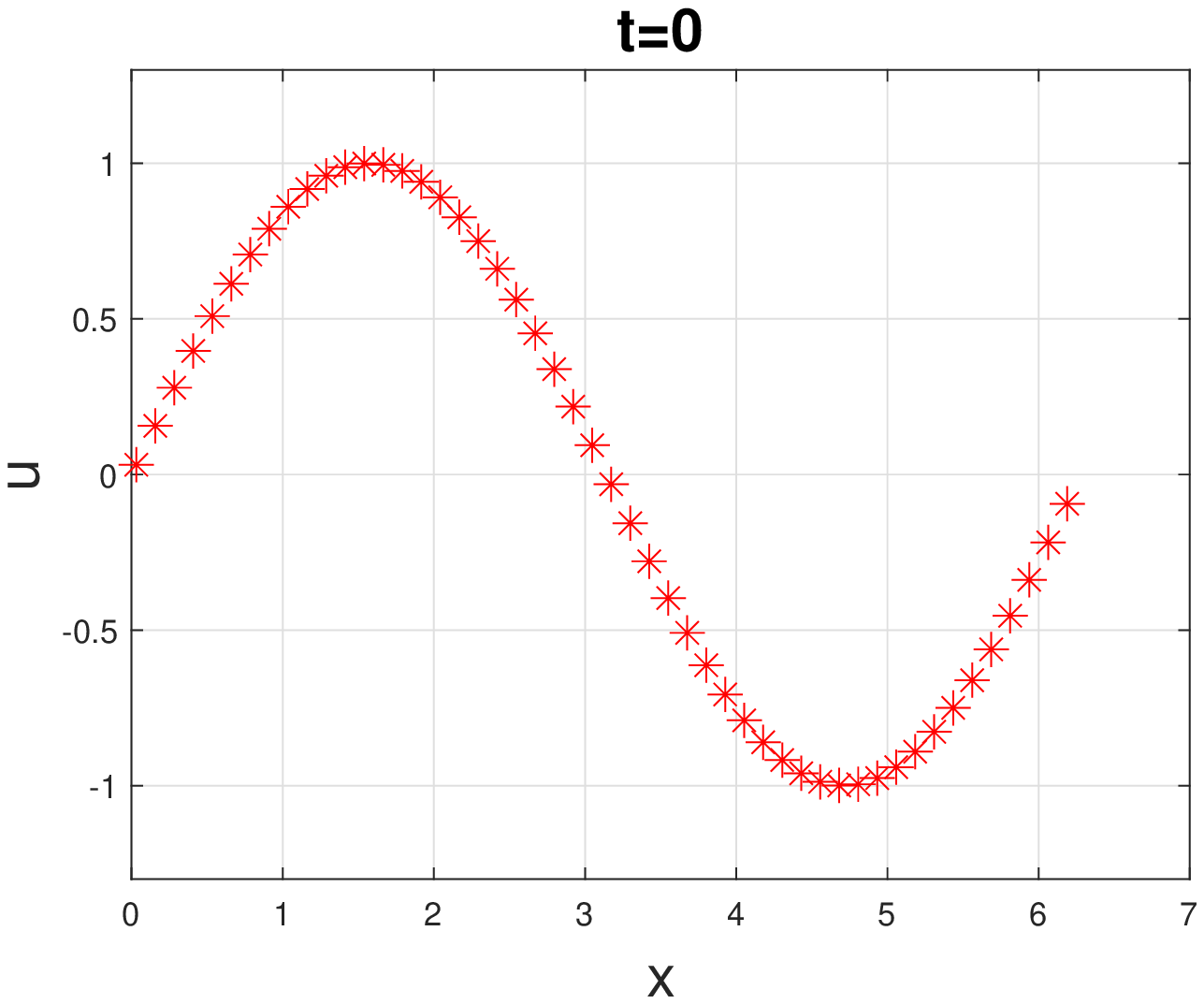}\quad\quad
\includegraphics[width=0.3\linewidth]{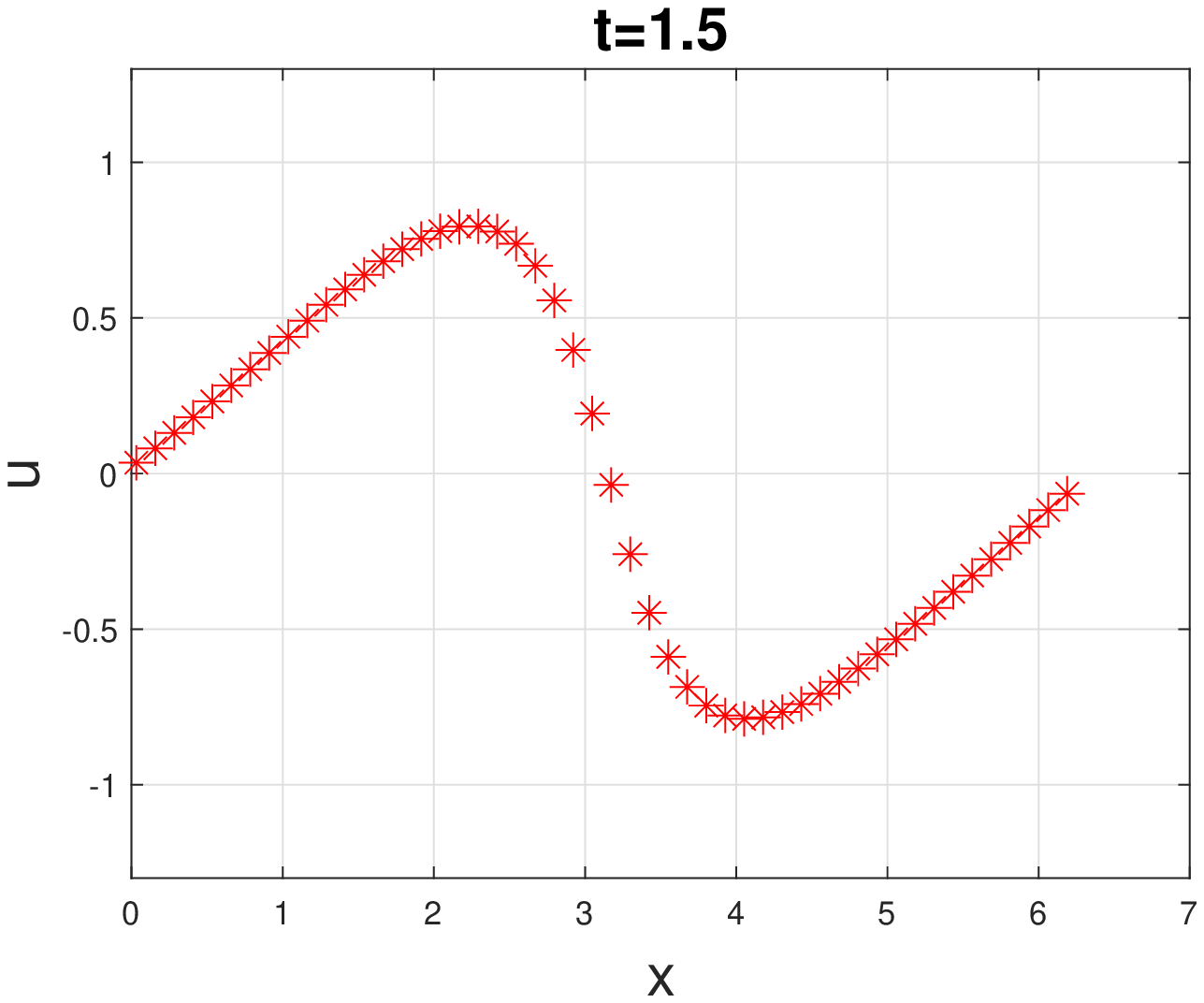}\\
\includegraphics[width=0.3\linewidth]{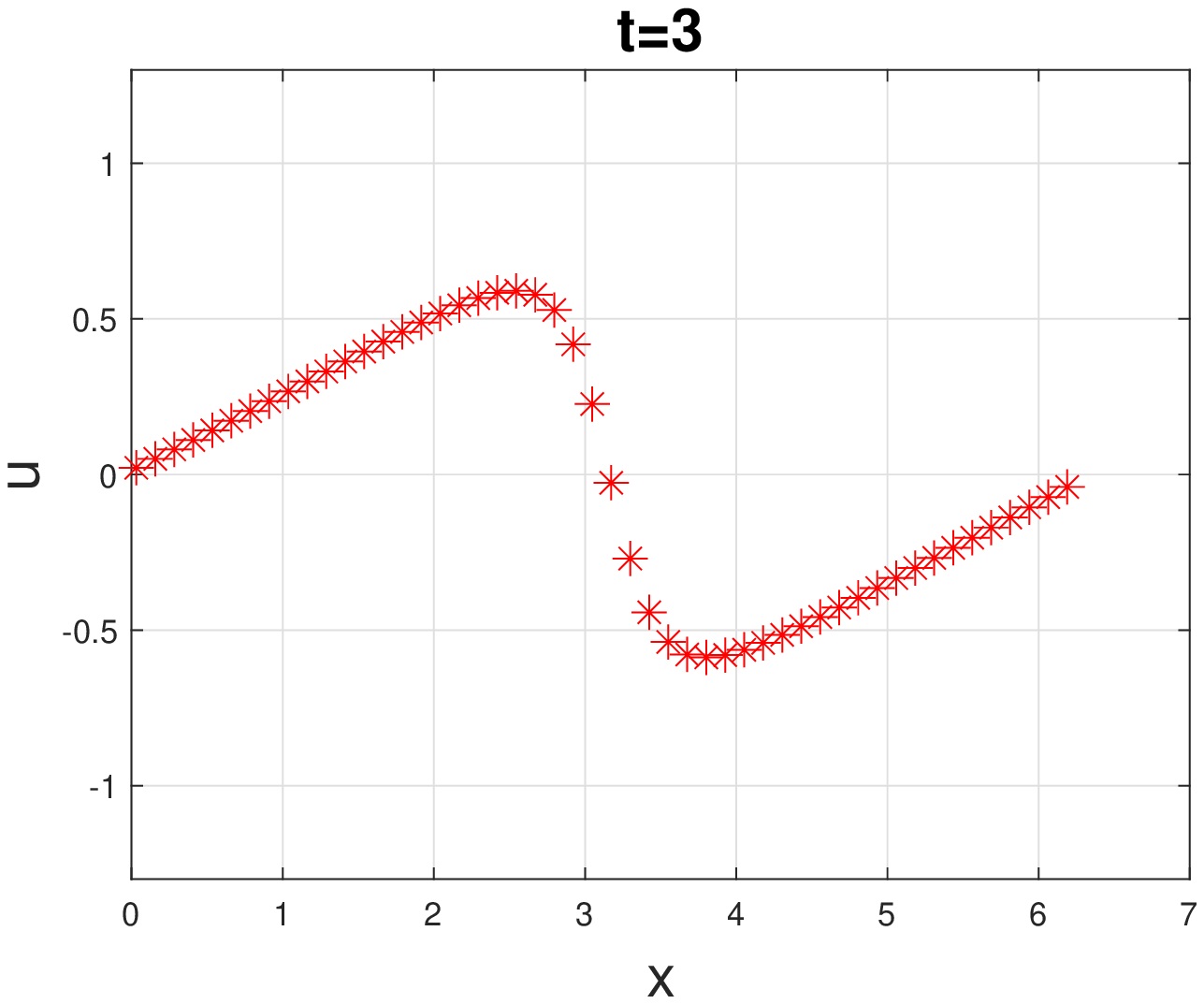}\quad\quad
\includegraphics[width=0.3\linewidth]{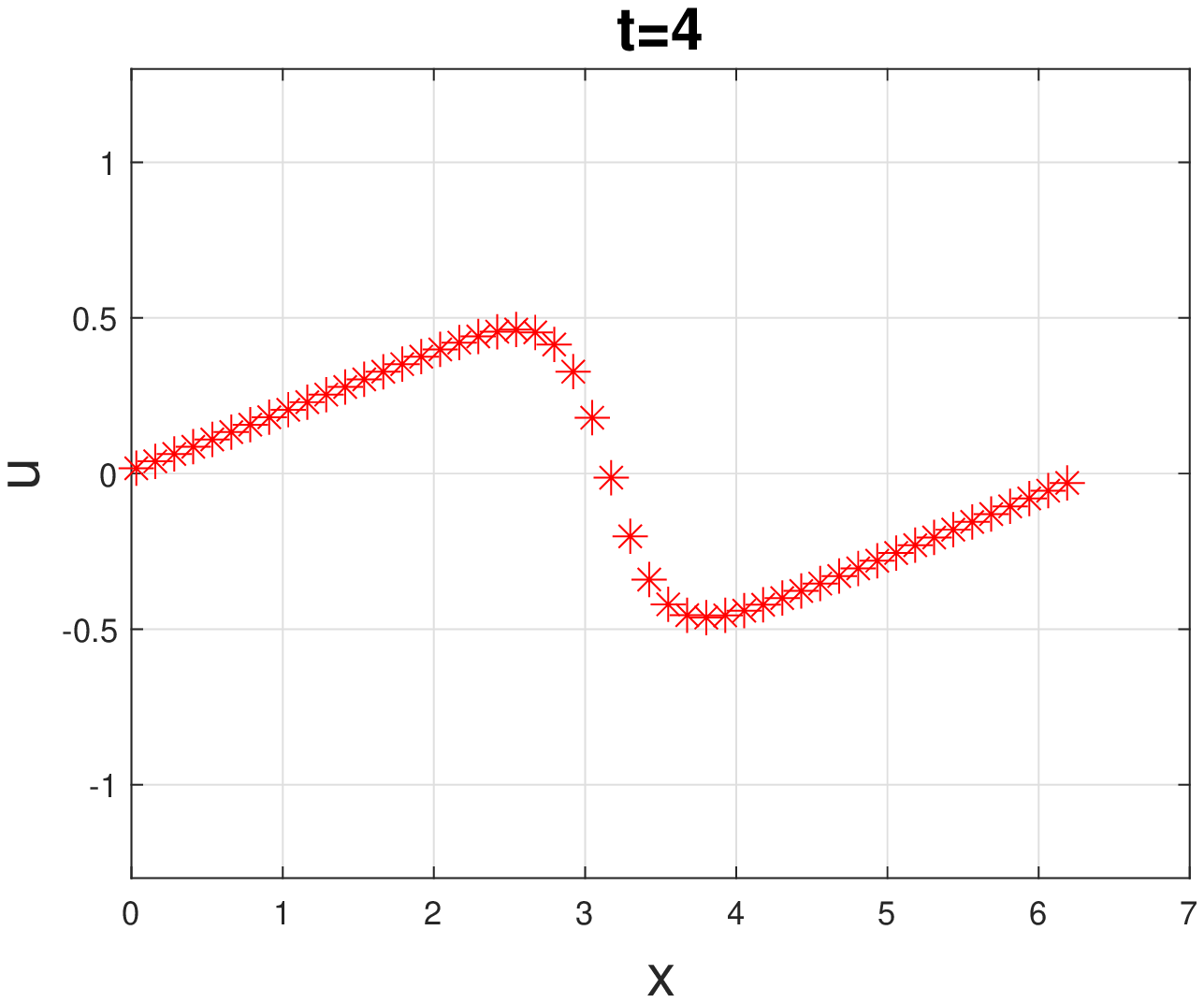}
\caption{viscous Burgers' equation (\ref{ex-viscous-burgers}) simulation with neural network method}
\label{fig:viscous-Burgers-sine}
\end{figure}

Solution evolution run by neural network method is shown in Figure \ref{fig:viscous-Burgers-sine}. The results are comparable to those obtained from classical numerical methods. 

\section{Conclusions}{\label{S:4}}

A finite volume or cell-average based neural network method is developed for hyperbolic and parabolic partial differential equations. Classical  numerical methods and principles motivate and guide the design of machine learning neural network methods. It is found the cell-average based neural network method can be relieved from the explicit scheme CFL restriction and is able to evolve the solution forward in time with large time step size. The method is able to accurately capture shock and rarefaction waves. We believe this is the right direction to explore neural networks solvers for partial differential equations, which is built upon solution properties and successful numerical methods.



%% References



%% New version of the num-names style

%\newpage

\bibliographystyle{elsarticle-num-names}
\bibliography{Neural_Network_for_PDE}






% \end{thebibliography}

\end{document}